\documentclass[12pt,twoside,a4paper]{amsart}
\UseRawInputEncoding
\usepackage{amsmath,amssymb}
\usepackage{cite}
\usepackage{graphicx}
\usepackage{txfonts}
\usepackage{mathrsfs}
\usepackage{empheq}
\newcommand{\ds}{\displaystyle}

\newcommand{\ka}{J.Kaczorowski}
\newcommand{\re}{M. Reko\'{s}}
\newcommand{\s}{\mathcal{S}^\text{poly}}
\newcommand{\se}{\mathcal{S}}
\newcommand{\f}{\overline{F(1-\overline{s})}}
\newcommand{\mf}{\overline{\mu_F (k)}}

\allowdisplaybreaks[4]
\theoremstyle{definition}
\newtheorem{definition}{Definition}[section]
\newtheorem{remark}[definition]{Remark}

\theoremstyle{theorem}
\newtheorem{theorem}{Theorem}[section]
\newtheorem{lemma}[theorem]{Lemma}

\newtheorem{fact}[theorem]{Fact}
\newtheorem{cor}[theorem]{Corollary}

\begin{document}
\begin{center}
\textbf{\large On some analytic properties of a function associated with the Selberg class satisfying certain special conditions}
\\
\vspace{0.9cm}
by
\\
\vspace{0.35cm}
Hideto  IWATA
\end{center}

\vspace{0.5cm}

\textbf{Abstract.} In 2001, {\re} described the analytic behavior for a function $f(z)$ connected with the Euler totient function for Im$\hspace{0.1cm}z > 0$ (see (1.2)) imitating the previous research of ~\cite{b} and ~\cite{k}. In the present paper, for Im$\hspace{0.1cm}z > 0$ we describe the analytic behavior of the generalized function $f(z,F)$ (see (2.1)), where the function $F$ belongs to the subclass of the Selberg class which has a polynomial Euler product and satisfies some special conditions.

\section{Introduction}
{\ka} defined the associated Euler totient function for a class of generalized $L$-functions including the Riemann zeta function, Dirichlet $L$-functions and obtained an asymptotic formula (see~\cite{Kac}):By a polynomial Euler product we mean a function $F(s)$ of a complex variable $s=\sigma+it$ which for $\sigma>1$ is defined by a product of the form
\begin{equation*}
   F(s) = \prod_{p} F_p (s) = \prod_{p}\prod_{j=1}^d \left( 1-\frac{\alpha_j (p)}{ p^s } \right)^{-1},   \tag{1.1}
\end{equation*}where $p$ runs over primes and $|\alpha_j (p)| \leq 1$ for all $p$ and $1 \leq j \leq d$. 
\footnote[0]
{\textit{2010 Mathematics Subject Classification. Primary}
 33C15, 11M06, 11M41 
   
\textit{Key words and phrases}
: polynomial Euler product,
Whittaker function,
Selberg class.
}We assume that $d$ is chosen as small as possible, i.e. that there exists at least one prime number $p_0$ such that
\[ \ds \prod_{j=1}^d \alpha_j ( p_0 ) \neq 0. \]Then $d$ is called the \textit{Euler degree} of $F$. Note that the $L$-functions from number theory including the Riemann zeta function, Dirichlet $L$-functions, Dedekind zeta and Hecke $L$-functions of algebraic number fields, as well as the (normalized) $L$-functions of holomorphic modular forms and, conjecturally, many other $L$-functions are polynomial Euler products. 

{\re} described the analytic property of some function connected with the Euler totient function (see~\cite{R}):We describe basic analytic properties of the function $f(z)$ defined for Im $z >0$ as follows :
\begin{equation*}
   f(z) = \lim_{n \to \infty} \sum_{\substack{\rho \\ 0< \text{Im} \> \rho < T_n}} \frac{e^{\rho z}\zeta(\rho-1)}{\zeta^\prime (\rho)},   \tag{1.2}
\end{equation*}where $T_n$ denotes a sequence of real numbers which yields appropriate grouping of the zeros. The summation is over non-trivial zeros of the Riemann zeta-function with positive imaginary part. For simplicity we assume here that the zeros are simple. {\re} showed the regularity of $f(z)$ for Im$\> z > 0$, meromorphic continuation to the whole complex plane, a certain functional equation and information of singularities. The functional equation for $f(z)$ connects the value of the function $f$ at the points $z$ and $\bar{z}$.  Let $\ell$ denote a smooth curve $\tau : [0,1] \longrightarrow \mathbb{C}$ such that $\tau(0) = -\frac{1}{4}, \> \tau(1) = \frac{5}{2}$ and $0<\text{Im} \> \tau <1$ for $t \in (0,1)$. The analytic property of $f(z)$ is described by the following theorems:
\begin{theorem}[Theorem 1. in ~\cite{R}]
The function $f(z)$ is analytic on the upper half-plane $\mathbb{H}$ and for $z \in \mathbb{H}$ we have 
\begin{equation*}
   2\pi if(z) = f_1 (z) + f_2 (z) - e^{\frac{5}{2}z} \sum_{n=1}^\infty \frac{\varphi(n)}{n^\frac{5}{2}(z-\log n)} \tag{1.3}
\end{equation*}where the last term on the right is a meromorphic function on the whole complex plane with the poles at $z=\log n, \> n=1,2,\ldots$. The function
\begin{equation*}
   f_1 (z) = \int_{-\frac{1}{4}+i\infty}^{-\frac{1}{4}} \frac{\zeta(s-1)}{\zeta(s)}e^{sz}ds   \tag{1.4}
\end{equation*}is analytic on $\mathbb{H}$ and
\begin{equation*}
   f_2 (z) = \int_{\ell(-\frac{1}{4},\frac{5}{2})} \frac{\zeta(s-1)}{\zeta(s)}e^{sz}ds   \tag{1.5}
\end{equation*}is analytic on the whole complex plane.
\end{theorem}

\begin{theorem}[Theorem 2. in ~\cite{R}]
The function $f(z)$ can be continued analytically to a meromorphic function on the whole complex plane, which satisfies the functional equation
\begin{equation*}
   f(z) + \overline{f(\bar{z})} = B(z)   \tag{1.6}
\end{equation*}and
\begin{equation*}
   B(z) = -\frac{6}{\pi^2}e^{2z} + \frac{1}{2\pi^2}\sum_{k,n=1}^\infty \frac{\mu(k)}{n^2 k}\left[ \frac{1}{(nke^z -1)^2} + \frac{2}{nke^z -1} + \frac{1}{(nke^z + 1)^2} -\frac{2}{nke^z +1} \right],   \tag{1.7}
\end{equation*}where $B(z)$ is a meromorphic function on the whole complex plane with the poles of the second order at $z=-\log nk, \> n,k=1,2,\ldots$ and $\mu(n)$ is the M\"{o}bius function.
\end{theorem}We see that the function $f(z)$ has  simple poles at $z=\log n \> (n=1,2,\ldots)$ with residue
\[ -\frac{\varphi(n)}{2\pi i}, \]where $\varphi(n)$ is the Euler totient function.

Now we provide the definition of the \textit{Selberg class} $\mathcal{S}$ as follows : $F \in \mathcal{S}$ if
\begin{enumerate}
\item[(i)](\textit{ordinary Dirichlet series}) $\displaystyle F(s) = \sum_{n=1}^\infty a_F (n)n^{-s}$, absolutely convergent for $\sigma > 1$; 
\item[(ii)](\textit{analytic continuation}) there exists an integer $m\geq0$ such that $(s-1)^m \cdot F(s)$ is an entire function of finite order; 
\item[(iii)](\textit{functional equation}) $F(s)$ satisfies a functional equation of type $\Phi(s) = \omega\overline{{\Phi}(1-\overline{s})}$, where
                                         \begin{equation*}
                                            \Phi(s) = Q^s \prod_{j=1}^r \Gamma(\lambda_j s + \mu_j)F(s) = \gamma(s)F(s),   \tag{1.8}
                                          \end{equation*}
                                            say, with $r\geq0, Q>0, \lambda_j >0$, Re\hspace{0.001cm} $\mu_j \geq 0$ and $|\omega| = 1$;
\item[(iv)](\textit{Ramanujan conjecture}) for every $\epsilon>0, a_F (n) \ll n^\epsilon$  
\item[(v)](\textit{Euler product}) $\displaystyle F(s) = \prod_{p} \exp\left( \sum_{\ell=0}^\infty \frac{b_F (p^\ell)}{p^{\ell s}} \right)$, where $b_F (n) = 0$ unless $n=p^m$ with $m\geq1$, and $b_F (n) \ll n^\vartheta$ for some 
                                  $\vartheta<\frac{1}{2}$.
\end{enumerate}Note that we understand an empty product is equal to $1$.

\section{The extension of $f(z)$ to the subclass of $\mathcal{S}$}
If a function $F \in \mathcal{S}$ has a polynomial Euler product (1.1), the subclass of $\mathcal{S}$ of the functions with polynomial Euler product is denoted by $\mathcal{S}^{\text{poly}}$. Obtaining results similar to Theorem 1.4.1 and Theorem 1.4.2 for a function $F$ belonging to $\mathcal{S}^{\text{poly}}$ is aim in this present paper. Now we assume that $(r, \lambda_j ) = (1,1)$ for all $j$ in the functional equation (1.8). The complex number $\mu_1$ when $r=1$ in (1.8) is hereafter referred to as $\mu$. Let $\rho$ denote the non-trivial zeros of $F$ with positive imaginary part. We assume that the order of $\rho$ is simple. Moreover, let $T_n$ denote a sequence of real numbers which yields appropriate grouping of the zeros which will be given later. For Im $z > 0$, we consider a function defined by
\begin{equation*}
   f(z,F) = \lim_{n \to \infty} \sum_{\substack{\rho \\ 0< \text{Im} \> \rho < T_n}} \frac{e^{\rho z}\zeta(\rho-1)}{F^\prime (\rho)}.   \tag{2.1}
\end{equation*}If there are trivial zeros of $F(s)$ on the imaginary axis in $\mathbb{H}$, we incorporate into the summation. The reason why the numerator on the right hand side is the Riemann-zeta function is because then we can use the Barnes type integral for the Whittaker function under the hypothesis $(r, \lambda_j ) = (1,1)$ for all $j$ in (1.8) (see Theorem 5.1 and Section 7). 
\begin{fact}
The series in (2.1) converges.
\end{fact}We will prove Fact 2.1 in the next section.
\section{Proof of Fact 2.1.}
We prove Fact 2.1. We use the following lemma.
\begin{lemma}[Lemma 4. in ~\cite{Sr}]
Let $F \in \mathcal{S}$ and let $T$ be sufficiently large. Moreover, let $H=D\log\log T$ be fixed, where $D$ is a large positive constant. In any subinterval of length $1$ in $[T-H, T+H]$ there are lines $t=t_0$ such that
\begin{equation*}
   |F(\sigma+it_0)|^{-1} = O(\exp(C(\log T)^2)),   \tag{3.1}
\end{equation*}uniformly in $\sigma\geq -2$, where $C$ is a positive constant.  
\end{lemma}Let $T$ be sufficiently large. We fix $H=D\log\log T$, where $D$ is a large positive constant. We take any subinterval $[n, n+1]$, where $n$  is a positive integer in $[T-H, T+H]$. Then, by Lemma 3.1 there are lines $t=T_n$ such that
\begin{align*}
|F(\sigma+iT_n)|^{-1} = O( \exp(C_1 (\log T)^2) )   \tag{3.2}
\end{align*}uniformly for $\sigma \geq -2$, where $C_1$ is a positive constant. Since $T_n$ is contained in the interval $[T-H,T+H]$, we can see that $T_n  \sim T$ as $n$ tends to infinity. Let $\alpha = \frac{1}{2}\min\{\text{Im}\hspace{0.1cm}\rho ; \text{Im}\hspace{0.1cm}\rho >0  \}$ and $\mathscr{L}$ denote the contour consisting of line segments
\[ \left[ b, b+i T_n \right], \> \left[ b+i T_n , a+i T_n \right], \> \left[ a+i T_n , a \right], \> \left[a,\frac{a+b}{2}+i\alpha \right], \> \left[\frac{a+b}{2}+i\alpha, b \right], \]where $\max \left\{ -\frac{3}{2}, \frac{1}{2}\max\left\{\text{Re} \hspace{0.1cm} \rho; \text{Re} \hspace{0.1cm} \rho < 0 \right\} \right\} < a < 0, b > \frac{5}{2}$. We assume that the real part of $s = a + it \> (t \in \mathbb{R})$ does not coincide the poles of $\Gamma( s + \mu )\Gamma(s - \mu)$. We consider the following contour integral round $\mathscr{L}$ :
\begin{equation*}
   \int_\mathscr{L} \frac{\zeta(s-1)}{F(s)}e^{sz}ds.  \tag{3.3}
\end{equation*} Since we assume the order of $\rho$ is simple, we have by residue theorem 
\begin{align*}
   \int_\mathscr{L} \frac{\zeta(s-1)}{F(s)}e^{sz}ds
   &= \int_{a+iT_n}^a \frac{\zeta(s-1)}{F(s)}e^{zs}ds + \int_L \frac{\zeta(s-1)}{F(s)}e^{zs}ds 
   \\
   &+ \int_b^{b+iT_n}\frac{\zeta(s-1)}{F(s)}e^{zs}ds + \int_{b+iT_n}^{a+iT_n}\frac{\zeta(s-1)}{F(s)}e^{zs}ds
   \\
   &= 2\pi i \sum_{\substack{\rho \\ 0< \text{Im} \> \rho < T_n}} \frac{e^{\rho z}\zeta(\rho-1)}{F^\prime (\rho)},   \tag{3.4}
\end{align*}where the path of integration $L$ consists of two line segments $\left[a,\frac{a+b}{2}+i\alpha \right]$ and $\left[\frac{a+b}{2}+i\alpha, b \right]$. We estimate the integral along the line segment $[b+iT_n , a+iT_n]$. By (3.2), for $a \leq \sigma \leq b$, we have the estimate
\[ |F(\sigma+iT_n)|^{-1} = O( \exp(C(\log T)^2) ). \] For $z=x+iy \> (y>0)$,
\begin{align*}
   &\left| \int_{a+iT_n}^{b+iT_n} \frac{\zeta(s-1)}{F(s)}e^{zs}ds \right|
   \\
   &\leq \int_{a+iT_n}^{b+iT_n} \left| \frac{\zeta(s-1)}{F(s)}e^{zs} \right| |ds|
   \\
   &\ll \int_{a}^b |\zeta(\sigma-1+iT_n)|\exp( C(\log T)^2 + x\sigma -yT_n )d\sigma
   \\
   &\ll (b-a)\exp\{ C(\log T)^2 -yT_n + |x|(|a| + |b|)\} T_{n}^c ,   \tag{3.5}
\end{align*}where the constant $c$ may depend on $a,b$. The last term on the right hand side in the above tends to zero as $n$ tends to infinity. By Theorem 4.1, the convergence of the other integrals in (3.4) are ensured (see (4.5)-(4.7)). Therefore, the series in (2.1) is convergent. $\square$ 

\section{Main theorems}
In (3.4), we have as $n$ tends to infinity 
\begin{equation*}
\int_{a+i\infty}^a \frac{\zeta(s-1)}{F(s)}e^{zs}ds + \int_L \frac{\zeta(s-1)}{F(s)}e^{zs}ds + \int_{b}^{b+i\infty} \frac{\zeta(s-1)}{F(s)}e^{zs}ds  = 2\pi if(z,F),   \tag{4.1}
\end{equation*}where the function $f(z,F)$ is defined in (2.1). To calculate the integral along the vertical line with $s=b+it \hspace{0.1cm} (t \geq 0)$, we prepare the Dirichlet series expansion of $\zeta(s-1)/F(s)$ for $\sigma>2$. 
\begin{definition}[p 34 in ~\cite{Kac}]
For $\sigma>1$ and $F \in {\s}$, we define the function $\mu_F$ as follows :
\begin{equation*}
   \frac{1}{F(s)} = \sum_{n=1}^\infty \frac{\mu_F (n)}{n^s} = \prod_p \prod_{j=1}^d \left( 1 - \frac{\alpha_j (p)}{p^s} \right).   \tag{4.2}
\end{equation*}
\end{definition}
\begin{remark}[p34 in ~\cite{Kac}]
By (4.2), $|\mu_F (n)| \leq \tau_d (n)$, where $\tau_d (n)$ is the divisor function of order $d$, so that $\zeta^d (s) = \sum_{n=1}^\infty \tau_d (n) / n^s$ for $\sigma>1$. In particular $\tau_1 (n) = 1$ for all $n$. 
\end{remark}Using (4.2), for $\sigma>2$
\begin{align*}
   \frac{\zeta(s-1)}{F(s)}
   &= \left( \sum_{l=1}^\infty \frac{\mu_F (l)}{l^s} \right)\left( \sum_{m=1}^\infty \frac{1}{m^{s-1}} \right) 
   \\
   &= \sum_{n=1}^\infty \frac{g(n)}{n^s},   \tag{4.3}
\end{align*}where
\begin{equation*}
   g(n) = \sum_{d|n} \mu_F (d)\frac{n}{d}.   \tag{4.4}
\end{equation*}
\begin{theorem}
The function (2.1) is analytic on $\mathbb{H}$ and for $z \in \mathbb{H}$ we have 
\begin{equation*}
   2\pi if(z,F) = f_1 (z,F) + f_2 (z,F) - e^{bz}\sum_{n=1}^\infty \frac{g(n)}{n^b (z-\log n)},   \tag{4.5}
\end{equation*}where the last term on the right is a meromorphic function on the whole complex plane with the poles at $z=\log n, \> n=1,2,\ldots$. The function
\begin{equation*}
   f_1 (z,F) = \int_{a+i\infty}^a \frac{\zeta(s-1)}{F(s)}e^{sz}ds   \tag{4.6}
\end{equation*}is analytic on $\mathbb{H}$ and 
\begin{equation*}
   f_2 (z,F) = \int_L \frac{\zeta(s-1)}{F(s)}e^{sz}ds   \tag{4.7}
\end{equation*} is analytic on the whole complex plane.
\end{theorem}
\begin{theorem}
For $F$ belonging to  ${\s}$ whose $(r,\lambda_j) = (1,1)$ for all $j$ in (1.8) and $0 \leq \mu < 1$, the function (2.1) has a meromorphic continuation to $y>-\pi$.
\end{theorem}
The $L$-functions associated with holomorphic cusp forms and Dedekind zeta functions of the imaginary quadratic fields are examples of $F$ considering in Theorem 4.2. Let
\begin{equation*}
   \mathbb{H}^{-} = \{ z \in \mathbb{C} : \text{Im}\hspace{0.1cm}  z < 0 \}.   \tag{4.8}
\end{equation*}We consider the function for $z \in \mathbb{H}^-$
\begin{equation*}
   f^{-} (z,F) = \lim_{n \to \infty} \sum_{\substack{\rho \\ -T_n < \text{Im} \> \rho < 0}} \frac{e^{\rho z}\zeta(\rho-1)}{F^\prime (\rho)}.   \tag{4.9}
\end{equation*}If there are trivial zeros of $F(s)$ on the imaginary axis in $\mathbb{H}^-$, we incorporate into the summation. The convergence for the series on the right hand side in (4.9) is proved by the same way as in Fact 2.1. 
\begin{cor}
For $F$ belonging to ${\s}$ which satisfies the same condition as in Theorem 4.2, the function (4.9) has a meromorphic continuation to $y<\pi$. 
\end{cor}
\begin{theorem}
For $F$ belonging to ${\s}$ which satisfies the same condition as in Theorem 4.2, the function (2.1) can be continued analytically on the whole complex plane. In addition to the condition as in Theorem 4.2, we assume that the coefficient $a_F (n)$ in the Dirichlet series is real value for all $n$. Then, the function (2.1) satisfies the functional equation
\begin{equation*}
   f(z,F) + \overline{f(\overline{z},F)} = B(z,F),  \tag{4.10}
\end{equation*}where
\begin{equation*}
   B(z,F) = \frac{1}{2\pi i}(f_1 (z,F) + f_{1}^{-}(z,F)) - \frac{e^{2z}}{F(2)}.   \tag{4.11}
\end{equation*}for all $z \in \mathbb{C}$. The definition of $f_{1}^- (z,F)$ is mentioned later (9.4). 
\end{theorem} If a function $F \in {\se}$ satisfies the conditions (i)-(iii) on ${\se}$, we denote the this class by ${\se}^\#$ and call the \textit{extended Selberg class}. It is known that the Dirichlet coefficient $a_F (n)$ of $F \in {\se}^{\#}$ which satisfies some special conditions is real (see~\cite{kp}).

\section{Some auxiliary results on the Whittaker function}
First, we recall the definition of the Whittaker function which is necessary to show the main theorems. To do this, we introduce some related functions. Secondly, we prepare some auxiliary results e.g.\hspace{-0.0005cm}the integral expression and the asymptotic expansion.
\begin{definition}[The confluent hypergeometric function (~\cite{L})]
Let $z$ be a complex variable, $\alpha$ and $\gamma$ are parameters which can take arbitrary real and complex values except that $\gamma \neq 0, -1, -2, \ldots$. Moreover,
\begin{equation*}
   (\lambda)_0 = 1, \quad (\lambda)_k = \frac{\Gamma(\lambda+k)}{\Gamma(\lambda)} = \lambda(\lambda+1)\cdots (\lambda+k-1) \quad (k \in \mathbb{N}).  \tag{5.1}
\end{equation*}We define the \textit{confluent hypergeometric function (Kummer's function)} as follows :
\begin{equation*} 
\Phi(\alpha, \gamma ; z) = \sum_{n=0}^\infty \frac{(\alpha)_n}{(\gamma)_n}\frac{z^n}{n!} \quad (|z|<\infty).  \tag{5.2}
\end{equation*} 
\end{definition}By ratio test, the series (5.2) is convergence absolutely for all $\alpha,\gamma$ and $z$, except $\gamma \neq 0,-1,-2,\ldots$. Hence, (5.2) is an analytic and one-valued function for all $z$. It is also to be noted that (5.2) is a particular solution of the linear differential equation ( \textit{Kummer's equation} )
\begin{equation*}
   z\frac{d^2 u}{dz^2} + (\gamma-z)\frac{du}{dz} - \alpha u = 0,   \tag{5.3}
\end{equation*}where $ \alpha, \gamma$ are the same as in (5.2).
\begin{definition}[The confluent hypergeometric function of the second kind (~\cite{L})]
We introduce a new function
\begin{gather*}
   \Psi(\alpha, \gamma ; z) = \frac{\Gamma(1-\gamma)}{\Gamma(1+\alpha-\gamma)}\Phi(\alpha,\gamma;z) + \frac{\Gamma(\gamma-1)}{\Gamma(\alpha)}z^{1-\gamma}\Phi(1+\alpha-\gamma,2-\gamma;z),
   \\
  (|\arg z| < \pi, \gamma \neq 0, \pm 1, \pm 2, \ldots)  \tag{5.4}
\end{gather*}called \textit{the confluent hypergeometric function of the second kind}. The condition $\gamma \neq 0, \pm 1, \pm 2, \ldots$ in (5.4) comes from the condition of the $\Gamma$-factors in the numerator, $\Phi(\alpha, \gamma ; z)$ and $\Phi(1+\alpha-\gamma,2-\gamma;z) $ on the right hand side in (5.4).
\end{definition}Since the function (5.4) is a many-valued function of $z$ for $\alpha$ and $\gamma$ real or complex, we take its branch which lies in the $z$-plane cut along the negative real axis. Also, (5.4) is analytic function for all $\alpha, \gamma$ and $z$ except $\gamma \neq 0,-1,-2,\ldots$.
\begin{definition}[The Whittaker function (~\cite{L})]
A class of functions related to the confluent hypergeometric functions, and often encountered in the applications, consists of the \textit{Whittaker function}, defined by the formula  
\begin{equation*}
   W_{k, \hspace{0.01cm} \mu}(z) = z^{\mu+\frac{1}{2}}e^{-\frac{z}{2}}\Psi\left( \frac{1}{2}-k+\mu, 2\mu+1 ; z \right) \quad (|\arg z| < \pi).  \tag{5.5}
\end{equation*}
\end{definition}By (5.4) or (5.5), it follows that (5.5) is a many-valued function of $z$. Therefore, we also take its branch as the same in (5.4).
\begin{theorem}[Barnes type integral for the Whittaker function (~\cite{Sl}, ~\cite{WW})]
The Barnes integral for $W_{k,\mu} (z)$ is
\begin{equation*}
W_{k,\mu}(z) = \frac{e^{-\frac{1}{2}z} z^k}{2\pi i}\int_{c-\infty i}^{c+\infty i}\frac{\Gamma(s)\Gamma\left( -s-k-\mu+\frac{1}{2} \right)\Gamma\left( -s-k+\mu+\frac{1}{2} \right)}{\Gamma\left( -k-\mu+\frac{1}{2} \right)\Gamma\left( -k+\mu+\frac{1}{2} \right)}z^s ds,   \tag{5.6}
\end{equation*}for $|\arg z| < \frac{3}{2}\pi$, and $k\pm \mu + \frac{1}{2} \neq 0,1,2,\ldots$; the contour has loops if necessary so that the poles of $\Gamma(s)$ and those of $\Gamma\left( -s-k-\mu +\frac{1}{2} \right) \Gamma\left( -s-k+\mu +\frac{1}{2} \right)$ are on opposite sides of it. 
\end{theorem}
In (5.6), it holds for all finite values of $c$ provided that the contour of integration can always be deformed so as to separate the poles $\Gamma(s)$ and those of the other $\Gamma$-factors. By Stirling's formula, the integral in (5.6) represents a function of $z$ which is analytic at all points in the domain $|\arg z| \leq \frac{3}{2}\pi - \alpha$, where $\alpha$ is any positive number.

\begin{theorem}[The asymptotic expansions in $z$ for $\Psi(a;b;z)$ (~\cite{Sl})]
We find that, as $z \to 0$,
\begin{align*}
   \Psi(a;b;z)
   &= \frac{\Gamma(b-1)}{\Gamma(a)}z^{1-b}+O(|z|^{\text{Re} b -2}) \quad (\text{Re} \hspace{0.1cm} b \geq 2, b \neq 2),   \tag{5.7}
   \\
   &= \frac{\Gamma(b-1)}{\Gamma(a)}z^{1-b}+O(|\log z|) \quad (b = 2),   \tag{5.8}
  \\
   &= \frac{\Gamma(b-1)}{\Gamma(a)}z^{1-b}+O(1) \quad (1 < \text{Re} \hspace{0.1cm} b < 2),   \tag{5.9}
   \\
   &= \frac{\Gamma(1-b)}{\Gamma(1+a-b)} + \frac{\Gamma(b-1)}{\Gamma(a)}z^{1-b} + O(|z|) \quad (\text{Re} \hspace{0.1cm} b =1, b \neq 1),   \tag{5.10}
   \\
   &= -\frac{1}{\Gamma(a)}\left\{ \log z + \frac{\Gamma^\prime}{\Gamma}(a) + 2C_0 \right\} + O(|z\log z|) \quad (b=1),   \tag{5.11} 
\end{align*}where $C_0$ is Euler's constant.
\end{theorem}
By the definition (5.5) and Theorem 5.2, we have the following asymptotic expansions in $z \to 0$ for $W_{k,\mu}(z)$. 
\begin{theorem}[The asymptotic expansions in $z \to 0$ for $W_{k,\mu} (z)$]
\begin{align*}
   W_{k,\mu}(z) 
   &= \frac{\Gamma(2\mu)}{\Gamma\left( \frac{1}{2} + \mu - k \right)}z^{\frac{1}{2}-\mu}+O(z^{\frac{3}{2}- \text{Re} \mu}) \quad \left( \text{Re} \hspace{0.1cm} \mu \geq \frac{1}{2}, \mu \neq \frac{1}{2} \right),  \tag{5.12}
   \\
   &= \frac{1}{\Gamma(1-k)} + O(|z\log z|) \quad \left( \mu = \frac{1}{2} \right),  \tag{5.13}
   \\
   &= \frac{\Gamma(2\mu)}{\Gamma\left( \frac{1}{2} + \mu - k \right)}z^{ \frac{1}{2} - \mu} + O(|z|^{\text{Re}\mu + \frac{1}{2}})  \quad \left( 0 < \text{Re} \hspace{0.1cm} \mu < \frac{1}{2} \right),  \tag{5.14}
   \\
   &= \frac{\Gamma(-2\mu)}{\Gamma\left( \frac{1}{2} - \mu -k \right)}z^{\mu+\frac{1}{2}} + \frac{\Gamma(2\mu)}{\Gamma\left( \mu+\frac{1}{2} - k \right)}z^{-\mu+\frac{1}{2}} + O(|z|^{\text{Re}\mu + \frac{3}{2}}) \quad 
   ( \text{Re} \hspace{0.1cm} \mu = 0, \mu \neq 0 ),  \tag{5.15}
   \\
   &= -\frac{z^\frac{1}{2}}{\Gamma\left( \frac{1}{2} - k \right)}\left( \log z + \frac{\Gamma^\prime}{\Gamma}\left( \frac{1}{2} - k \right) + 2C_0 \right) + O(|z|^\frac{3}{2} |\log z|) \quad (\mu = 0).   \tag{5.16}
\end{align*}
\end{theorem}
\section{Proof of Theorem 4.1.}
By (4.1), for $z \in \mathbb{H}$ we have
\[ f_1 (z,F) + f_2 (z,F) + f_3 (z,F) = 2\pi if(z,F), \] where $f_1 (z,F),  f_2 (z,F),  f_3 (z,F)$ denote corresponding integrals in (4.1), respectively. First, we calculate the integral
\[ f_3 (z,F) = \int_{b}^{b+i\infty} \frac{\zeta(s-1)}{F(s)}e^{zs}ds. \] Using the Dirichlet series expansion (4.3), we have
\begin{align*}
   f_3 (z,F)
   &= \int_{b}^{b+i\infty} \left( \sum_{n=1}^\infty \frac{g(n)}{n^s} \right)e^{zs}ds
   \\
   &= \sum_{n=1}^\infty g(n) \int_{b}^{b+i\infty}e^{s(z-\log n)}ds
   \\
   &= -e^{bz}\sum_{n=1}^\infty \frac{g(n)}{n^b (z-\log n)}.   \tag{6.1}
\end{align*}The interchange of the order of iintegration and summation is justified for $z \in \mathbb{H}$ by the absolute and uniform convergence of the series on the third line in (6.1). 

Secondly, we prove the function (4.7) is analytic on the whole complex plane. Since the length of $L$ is finite, $\zeta(s-1) = O(1)$. Also, there are no zeros on $L$, the function $\{F(s)\}^{-1}$ is bounded. Therefore, for $z=x+iy$ and $s=\sigma+it \> (a \leq \sigma \leq b, 0 \leq t \leq \alpha)$,
\begin{align*}
   \int_L \left| \frac{\zeta(s-1)}{F(s)}e^{zs} \right| |ds|
   &\ll \int_L \exp(x\sigma-yt)|\zeta(s-1)||ds|
   \\
   &= O\left( \int_L \exp(x\sigma-yt)|ds| \right).
\end{align*}In the case $y\geq 0$,
    \[ \int_L \exp(x\sigma-yt)|ds| \ll \int_L \exp(x\sigma)|ds| \ll_x 1. \]In the case $y<0$, 
    \[ \int_L \exp(x\sigma-yt)|ds| \ll_{x,y} 1. \]Hence the function (4.7) is analytic on the whole complex plane. 
    
Finally, we prove the function (4.6) is analytic on $\mathbb{H}$. Here we use the following two lemmas related to the Selberg class.  
\begin{definition}[p29 in ~\cite{P}]
Let $F \in \mathcal{S}$, 
\begin{equation*}d_F = 2\sum_{j=1}^r \lambda_j   \tag{6.2}
\end{equation*} is the \textit{degree} of $F(s)$. 
\end{definition}
\begin{lemma}[Theorem 3.1 in ~\cite{C and G}]
If $F \in \mathcal{S}$, then $F=1$ or $d_F \geq 1$.
\end{lemma}
\begin{lemma}[~\cite{M}.(8), p423]
For $\sigma<0$,
\begin{equation*}
  |F(\sigma+it)| \asymp t^{d_F (\frac{1}{2}-\sigma)}|F(1-\sigma+it)|   \tag{6.3}
\end{equation*}as $t \to \infty$.
\end{lemma}Using Lemma 6.2, we have
\begin{equation*}
     t^{\left( \frac{1}{2}-a \right)d_F} |F(1-a+it)| \ll |F(a+it)| \ll t^{\left( \frac{1}{2}-a \right)d_F}|F(1-a+it)|.
\end{equation*}Since $\max \left\{ -\frac{3}{2}, \frac{1}{2}\max\{ \text{Re} \hspace{0.1cm} \rho ; \text{Re} \hspace{0.1cm} \rho < 0 \} \right\} < a < 0$, so $\{F(1-a+it)\}^{-1}$ is bounded. Hence, we have
\[ |F(a+it)|^{-1} \ll t^{-(\frac{1}{2}-a)d_F}. \]
Also, by the functional equation for $\zeta(s)$,
\begin{align*}
   \zeta(s-1) 
   &= -\frac{1}{\pi}(2\pi)^{s-1}\cos\left( \frac{\pi}{2}s \right)\Gamma(2-s)\zeta(2-s)
   \\
   &\ll |t|^{\frac{3}{2}-a}.
\end{align*}In the case $F=1$, we have
\begin{align*}
   \int_{a+i\infty}^a \left| \frac{\zeta(s-1)}{F(s)}e^{zs} \right||ds|
   &\ll \int_{0}^\infty t^{\frac{3}{2}-a}e^{ax-yt}dt
   \\
   &\ll_{x,y} 1
   \\
   &\ll 1.
\end{align*}Next consider the case $F \neq 1$. Now, $F(s)$ belongs to $\mathcal{S}^{\text{poly}}$, so $F(s)$ is analytic except the point $s=1$. Since $\max \left\{ -\frac{3}{2}, \frac{1}{2}\max\{ \text{Re} \hspace{0.1cm} \rho ; \text{Re} \hspace{0.1cm} \rho < 0 \} \right\} < a < 0$, $|F(a+it)|^{-1}$ is bounded near $t=0$. Also, $d_F \geq 1$ by Lemma 6.1, we have 
\begin{align*}
   \int_{a+i\infty}^a \left| \frac{\zeta(s-1)}{F(s)}e^{zs} \right||ds|
   &\ll \int_{0}^\infty \frac{e^{ax-yt}}{|F(a+it)|}t^{\frac{3}{2}-a}dt
   \\
   &\ll \int_{0}^\infty t^{-\left( \frac{1}{2} - a \right)d_F + \frac{3}{2} - a }e^{ax-yt}dt
   \\
   &< \int_{0}^\infty t^{-\frac{1}{2}d_F + \frac{3}{2} - a }e^{ax-yt}dt.
\end{align*}Therefore, the integral on the above third line is absolutely and uniformly convergent on every compact subset of $\mathbb{H}$. Consequently, the function (4.6) is analytic for $y =$ Im $z >0$. $\square$
\section{Proof of Theorem 4.2}
We prove that the function $f(z,F)\>(z=x+iy)$ has a meromorphic continuation to $y>-\pi$. By Theorem 4.1, the function
\begin{align*} 
   f_1 (z,F) 
   &= \int_{a+i\infty}^a \frac{\zeta(s-1)}{F(s)}e^{zs}ds
   \\
   &= -\int_a^{a + i\infty} \frac{\zeta(s-1)}{F(s)}e^{zs}ds
\end{align*}is convergent for $y>0$. We recall the hypotheses that $(r, \lambda_j ) = (1,1)$ for all $j$ in (1.8) and $0 \leq \mu < 1$. We rewrite the functional equation (1.8) under these hypotheses as follows :
\begin{align*}
   Q^s \Gamma( s + \mu )F(s) 
   &= \omega \overline{Q^{1-\overline{s}}\Gamma( 1 - \overline{s} + \mu )F(1-\overline{s})}
   \\
   &= \omega Q^{1-s} \Gamma( 1 - s + \mu ){\f},
\end{align*}where the conditions of $Q$ and $\omega$ are the same as noted in (1.8). Hence
\begin{equation*} 
   \frac{1}{F(s)} 
   = \overline{\omega}Q^{2s-1}\frac{\Gamma( s + \mu )}{\Gamma( 1- s + \mu )}\frac{1}{{\f}}.   \tag{7.1}
\end{equation*}Using the following elementary formula for the $\Gamma$-function 
\[ \Gamma(s)\Gamma(1-s) = \frac{\pi}{\sin \pi s}, \](7.1) yields 
\begin{equation*}
   \frac{1}{F(s)} = \frac{\overline{\omega}}{\pi}Q^{2s-1}\sin\pi(s-\mu)\Gamma(s+\mu)\Gamma(s-\mu)\frac{1}{{\f}}.  \tag{7.2}
\end{equation*}By (7.2) and the functional equation for $\zeta(s)$, we have
\begin{align*}
   &f_1 (z,F)
   \\ 
   &=  -\int_a^{a + i\infty} \frac{\zeta(s-1)}{F(s)}e^{zs}ds
   \\
   &= \frac{\overline{\omega}}{2\pi^3 Q}\int_a^{a+i\infty} \hspace{-0.7cm} (2\pi Q^2)^s \cos\left( \frac{s}{2}\pi \right) \Gamma(2-s)\zeta(2-s)\sin\pi(s-\mu)
   \\
   &\times \Gamma(s+\mu)\Gamma(s-\mu)\frac{e^{zs}}{{\f}}ds
   \\
   &= \frac{\overline{\omega} e^{-\mu\pi i}}{(2\pi)^3 Qi}\int_a^{a+i\infty} (2\pi Q^2 )^s \frac{\zeta(2-s)}{{\f}}\Gamma(s - \mu)\Gamma(s + \mu )\Gamma(2-s)e^{\left( z + \frac{3}{2}\pi i \right)s}ds
   \\
   &-\frac{\overline{\omega} e^{\mu\pi i}}{(2\pi)^3 Qi}\int_a^{a+i\infty} (2\pi Q^2 )^s \frac{\zeta(2-s)}{{\f}}\Gamma( s - \mu )\Gamma(s + \mu)\Gamma(2-s)e^{\left(z - \frac{\pi}{2}i \right)s}ds   
   \\
   &+\frac{\overline{\omega} e^{-\mu\pi i} }{(2\pi)^3 Qi}\int_a^{a+i\infty} (2\pi Q^2 )^s \frac{\zeta(2-s)}{{\f}}\Gamma( s - \mu )\Gamma(s + \mu)\Gamma(2-s)e^{\left(z + \frac{\pi}{2}i \right)s}ds   
   \\   
   &-\frac{\overline{\omega} e^{\mu\pi i} }{(2\pi)^3 Qi} \int_a^{a+i\infty} (2\pi Q^2 )^s \frac{\zeta(2-s)}{{\f}}\Gamma(s - \mu )\Gamma( s + \mu )\Gamma(2-s) e^{(z - \frac{3}{2}\pi i)s}ds   
   \\
   &= f_{11}(z,F) + f_{12}(z,F) + f_{13}(z,F) + f_{14}(z,F),  \tag{7.3}
\end{align*}where $f_{11}(z,F), f_{12}(z,F), f_{13}(z,F), f_{14}(z,F)$ denote the corresponding integrals in (7.3) respectively. By Stirling's formula
\begin{equation*}
   \Gamma(s+\mu)\Gamma(s-\mu)\Gamma(2-s) \asymp e^{-\frac{3}{2}\pi |t|}|t|^{a+\frac{1}{2}}
\end{equation*}as $t$ tends to infinity and $\max \left\{ -\frac{3}{2}, \frac{1}{2}\max\{ \text{Re} \hspace{0.1cm} \rho ; \text{Re} \hspace{0.1cm} \rho < 0\} \right\} < a < 0$, $f_{11}(z,F)$ is analytic for $y > -3\pi$, $f_{12}(z,F)$ for $y> -\pi$, $f_{13}(z,F)$ for $y > -2\pi$, $f_{14}(z,F)$ for $y>0$. Splitting the integral in $f_{14}(z,F)$, we have
\begin{align*}
   f_{14}(z,F)
   &= -\frac{\overline{\omega} e^{\mu\pi i} }{(2\pi)^3 Qi} \int_a^{a+i\infty} \hspace{-0.5cm} (2\pi Q^2 )^s \frac{\zeta(2-s)}{{\f}}\Gamma(s - \mu )\Gamma( s + \mu )\Gamma(2-s) e^{(z - \frac{3}{2}\pi i)s}ds   
   \\   
   &= \frac{\overline{\omega}e^{\mu\pi i} i}{(2\pi)^3 Q}\left\{ \int_{a-i\infty}^{a+i\infty} - \int_{a-i\infty}^a \right\} (2\pi Q^2 )^s \frac{\zeta(2-s)}{{\f}}
   \\
   &\times \Gamma(s - \mu )\Gamma( s + \mu )\Gamma(2-s) e^{(z - \frac{3}{2}\pi i)s}ds.
\end{align*}We consider
\begin{equation*}
   I_1 (z,F) = \int_{a-i\infty}^{a+i\infty} (2\pi Q^2 )^s \frac{\zeta(2-s)}{{\f}}\Gamma(s - \mu )\Gamma( s + \mu )\Gamma(2-s) e^{(z - \frac{3}{2}\pi i)s}ds   \tag{7.4}
\end{equation*}and
\begin{equation*}
  I_2 (z,F) =  \int_{a-i\infty}^a (2\pi Q^2 )^s \frac{\zeta(2-s)}{{\f}}\Gamma(s - \mu )\Gamma( s + \mu )\Gamma(2-s) e^{(z - \frac{3}{2}\pi i)s}ds.   \tag{7.5}
\end{equation*}Here, we recall the hypothesis that the real part of $s=a+it \> (t \in \mathbb{R})$ does not coincide with the poles of $\Gamma(s+\mu)\Gamma(s-\mu)$. By Stirling's formula,  we can see that the integral $I_2 (z,F)$ is convergent for $y<3\pi$. Since the function $f_{14}(z,F)$ was analytic for $y>0$, the integral $I_1 (z,F)$ is convergent for $0<y<3\pi$. Using the Dirichlet series expansion (4.2), we have formally
\begin{align*}
   I_1 (z,F)
   &=  \int_{a-i\infty}^{a+i\infty} \hspace{-0.5cm} (2\pi Q^2)^s \left( \sum_{n=1}^\infty \frac{1}{n^{2-s}} \right) \overline{\left( \sum_{k=1}^\infty \frac{\mu_F (k)}{k^{1-\overline{s}}} \right)} \Gamma( s - \mu )\Gamma( s + \mu )\Gamma(2-s) e^{(z - \frac{3}{2}\pi i)s}ds
   \\
   &= \sum_{k,n=1}^\infty \frac{{\mf}}{kn^2} \int_{a-i\infty}^{a+i\infty} e^{s\left\{ \log (2\pi nk Q^2) -\frac{3}{2}\pi i + z \right\}}\Gamma( s - \mu )\Gamma( s + \mu )\Gamma(2-s)ds.  \tag{7.6}
\end{align*}The justification of the interchange of the order of integration and summation is ensured as follows : For $0<y<3\pi$, we have 
\begin{eqnarray*}
   \lefteqn{\int_{a-i\infty}^{a+i\infty} \left| e^{s\left\{ \log (2nk\pi Q^2) + z - \frac{3}{2}\pi i \right\}}  \Gamma( s - \mu )\Gamma( s + \mu )\Gamma(2-s)  \right| |ds|}
   \\
   &&\ll (nk)^a \cdot e^{ax} \int_{-\infty}^\infty e^{-\left( y - \frac{3}{2}\pi \right)t -\frac{3}{2}\pi|t|} |t|^{a+\frac{1}{2}}dt
   \\
   &&= (nk)^a \cdot e^{ax} \left\{ \int_{-\infty}^0 e^{-\left( y - \frac{3}{2}\pi \right)t + \frac{3}{2}\pi t} (-t)^{a+\frac{1}{2}}dt + \int_{0}^\infty e^{-\left( y - \frac{3}{2}\pi \right)t - \frac{3}{2}\pi t}t^{a+\frac{1}{2}}dt  \right\}
   \\
   &&\ll_{a,x,y} (nk)^a
   \\
   &&\ll_a  (nk)^a .
\end{eqnarray*}Hence, we have
\begin{equation*}
   \sum_{k,n=1}^\infty \left| \frac{{\mf}}{kn^2} \int_{a-i\infty}^{a+i\infty} e^{s\left\{ \log (2nk\pi Q^2) + z - \frac{3}{2}\pi i \right\}}\Gamma( s - \mu )\Gamma( s + \mu )\Gamma(2-s)ds \right|
   \ll_a  \sum_{k,n=1}^\infty \left| \frac{{\mf}}{k^{1-a} n^{2-a}} \right|.
\end{equation*}Since the series for $k$ is absolutely convergent by (4.2), the above series of the left hand side is convergent absolutely and uniformly. Therefore, the interchange of the order of integration and summation is justified for $0 < y < 3\pi$. Let $m_1 ,m_2$ be non-negative integers. The residue of the integrand in (7.6) at $s = \mu - m_1$ is 
\begin{align*}
   R^{(1)}_{k,n, m_1}(z,\mu)
   &= \lim_{s \to \mu - m_1} \{ s - ( \mu - m_1 ) \}\Gamma(s-\mu)\Gamma(s+\mu)\Gamma(2-s)e^{ \left\{ \log(2\pi nkQ^2)-\frac{3}{2}\pi i + z \right\}s }
   \\
   &= \frac{(-1)^{m_1}}{m_{1}!}\Gamma(2\mu-m_1)\Gamma(2-\mu+m_1)(2\pi nkQ^2)^{\mu-m_1}e^{\left( z - \frac{3}{2}\pi i \right)( \mu - m_1 )}.  \tag{7.7}
\end{align*}Similarly, the residue of the integrand in (7.6) at $s = -\mu - m_2$ is 
\begin{equation*}
   R^{(1)}_{k,n, m_2}(z,\mu)
   = \frac{(-1)^{m_2}}{m_{2}!}\Gamma(-2\mu-m_2)\Gamma(2+\mu+m_2)(2\pi nkQ^2)^{-\mu-m_2}e^{\left( z - \frac{3}{2}\pi i \right)( -\mu - m_2 )}.  \tag{7.8}
\end{equation*}When $\mu=0$, $\{\Gamma(s)\}^2$ has a double pole at $s=-m$, where $m$ is a non-negative integer. For every positive $\epsilon$, using the Taylor expansion for the every factor of the integrand in (7.6) at $s=-m+\epsilon$, the residue of the integrand in (7.6) at $s=-m$ is
\begin{equation*}
   R^{(1)}_{k,n,m,0}(z)
   = \frac{m+1}{m!}\frac{e^{-\left( z-\frac{3}{2}\pi i \right)m}}{(2\pi nkQ^2)^m}\left\{ \log(2\pi nkQ^2) + z -\frac{3}{2}\pi i + \sum_{k_1 = 1}^m \frac{1}{k_1} -C_0 - \frac{1}{m+1} \right\}.   \tag{7.9}
\end{equation*}Similarly, when $\mu=\frac{1}{2}$, $\Gamma\left( s-\frac{1}{2} \right)\Gamma\left( s+\frac{1}{2} \right) = \left( s-\frac{1}{2} \right)\left\{ \Gamma\left( s-\frac{1}{2} \right) \right\}^2$ has a double pole at $s=\frac{1}{2}-m$. In the same way, the residue of the integrand in (7.6) at $s=\frac{1}{2}-m$ is
\begin{align*}
   R^{(1)}_{k,n,m,\frac{1}{2}}(z)
   &= \frac{\Gamma\left( \frac{3}{2}+m \right)}{(m!)^2}\frac{e^{\left( \frac{1}{2}-m \right)\left( z -\frac{3}{2}\pi i \right)}}{(2\pi nkQ^2)^{m-\frac{1}{2}}}\left\{ m\left( \psi\left( \frac{3}{2}+m \right) -2\psi(m+1) \right)  \right.
   \\
   &- \left. m\left( \log (2\pi nkQ^2) + z -\frac{3}{2}\pi i \right) + 1 \right\},   \tag{7.10}
\end{align*}where $\psi(s)$ is the logarithmic derivative of $\Gamma(s)$, i.e. 
\begin{equation*}
   \psi(s) := \frac{\Gamma^\prime}{\Gamma}(s).  \tag{7.11}
\end{equation*}Since for $0<y<3\pi$ integrals along the upper and the lower side of the contour tend to $0$, for $\mu$ except $\mu \neq 0,\frac{1}{2}$, we have by theorem of residues
\footnotesize
\begin{align*}
   &\int_{a-i\infty}^{a+i\infty}e^{ s\left\{ \log (2nk\pi Q^2) + z - \frac{3}{2}\pi i \right\} } \Gamma(s - \mu )\Gamma( s + \mu )\Gamma(2-s)ds
   \\
   &= -\int_C e^{ s\left\{ \log (2nk\pi Q^2) + z - \frac{3}{2}\pi i \right\} } \Gamma(s - \mu )\Gamma( s + \mu )\Gamma(2-s)ds
   \\
   &-2\pi i\left\{ \sum_{m_1 = 0}^M R^{(1)}_{k,n, m_1}(z,\mu) + \sum_{m_2 = 0}^{M^\prime} R^{(1)}_{k,n, m_2}(z,\mu)  \right\},   \tag{7.12}
\end{align*}
\normalsize
where $C$ is the contour which the poles of $\Gamma(2-s)$ and those of $\Gamma( s + \mu) \Gamma ( s - \mu )$ are on opposite sides of it. Of course, when $\mu=0,\frac{1}{2}$, the terms of residue in (7.12) are replaced by
\begin{equation*}
   \sum_{m = 0}^{M^{\prime\prime}} R^{(1)}_{k,n,m,0}(z), \quad \sum_{m = 0}^{M^{\prime\prime}} R^{(1)}_{k,n,m,\frac{1}{2}}(z).   \tag{7.13} 
\end{equation*}We use the same convention hereafter. Putting $w=2-s$ in the integral round the contour $C$ on the right hand side in (7.12) and using Barnes type integral for the Whittaker function (5.6), for $\mu$ which is not integer except zero and $|y-\frac{3}{2}\pi|<\frac{3}{2}\pi$, the integral round the contour $C$ in (7.12) is
\begin{equation*}
   2\pi i(2\pi nkQ^2)^\frac{1}{2}\exp\left( \frac{e^{\frac{3}{2}\pi i - z}}{4\pi nkQ^2} + \frac{z}{2} - \frac{3}{4}\pi i \right)\Gamma(2-\mu)\Gamma(2+\mu)
      W_{-\frac{3}{2}, \mu}\left( \frac{e^{\frac{3}{2}\pi i - z}}{2\pi nkQ^2} \right).  \tag{7.14}
\end{equation*}Therefore, we have
\footnotesize
\begin{align*}
   I_1 (z,F)
  &= \sum_{k,n=1}^\infty \frac{{\mf}}{kn^2} \left\{ 2\pi i \cdot (2\pi nkQ^2)^\frac{1}{2}\exp\left( \frac{e^{\frac{3}{2}\pi i - z}}{4\pi nkQ^2} + \frac{z}{2} - \frac{3}{4}\pi i \right) \right.
  \\
  &\times \left. \Gamma(2-\mu)\Gamma(2+\mu)W_{-\frac{3}{2}, \mu}\left( \frac{e^{\frac{3}{2}\pi i - z}}{2\pi nkQ^2} \right) -2\pi i \left(\sum_{m_1 = 0}^M R^{(1)}_{k,n, m_1}(z,\mu) + \sum_{m_2 = 0}^{M^\prime} R^{(1)}_{k,n, m_2}(z,\mu) \right) \right\}.   \tag{7.15}
\end{align*}
\normalsize
The following lemma ensures the convergence for the series on the right hand side in (7.15).
\begin{lemma}
The series on the right hand side in (7.15) is absolutely and uniformly convergent on every compact subset on the whole complex plane.
\end{lemma}
We will prove Lemma 7.1 in the next section. By Lemma 7.1, for $F \in \mathcal{S}^{\text{poly}}$ whose $(r,\lambda_j) = (1,1)$ for all $j$ in (1.8) and $0 \leq \mu  < 1,$ we have the following analytic continuation of $f_1 (z,F)$ for $y > -\pi$ :
\footnotesize
\begin{align*}
   f_1 (z,F) 
   &= \frac{\overline{\omega} e^{-\mu\pi i}}{(2\pi)^3 Qi}\int_a^{a+i\infty} (2\pi Q^2 )^s \frac{\zeta(2-s)}{{\f}}\Gamma(s - \mu)\Gamma(s + \mu )\Gamma(2-s)e^{\left( z + \frac{3}{2}\pi i \right)s}ds
   \\
   &-\frac{\overline{\omega} e^{\mu\pi i}}{(2\pi)^3 Qi}\int_a^{a+i\infty} (2\pi Q^2 )^s \frac{\zeta(2-s)}{{\f}}\Gamma( s - \mu )\Gamma(s + \mu)\Gamma(2-s)e^{\left(z - \frac{\pi}{2}i \right)s}ds   
   \\
   &+\frac{\overline{\omega} e^{-\mu\pi i} }{(2\pi)^3 Qi}\int_a^{a+i\infty} (2\pi Q^2 )^s \frac{\zeta(2-s)}{{\f}}\Gamma( s - \mu )\Gamma(s + \mu)\Gamma(2-s)e^{\left(z + \frac{\pi}{2}i \right)s}ds   
   \\
   &+ \frac{\overline{\omega} e^{\mu\pi i}i }{(2\pi)^3 Q} \sum_{k,n=1}^\infty \frac{{\mf}}{kn^2}\times 2\pi i\left\{(2\pi nkQ^2)^\frac{1}{2}\exp\left( \frac{e^{\frac{3}{2}\pi i - z}}{4\pi nkQ^2} + \frac{z}{2} - \frac{3}{4}\pi i \right) \right.
   \\
   &\times \left. \Gamma(2-\mu)\Gamma(2+\mu)W_{-\frac{3}{2}, \mu}\left( \frac{e^{\frac{3}{2}\pi i - z}}{2\pi nkQ^2} \right) -\sum_{m_1 = 0}^M R^{(1)}_{k,n, m_1}(z,\mu) - \sum_{m_2 = 0}^{M^\prime} R^{(1)}_{k,n, m_2}(z,\mu) \right\}
   \\
   &+ \frac{\overline{\omega} e^{\mu\pi i} }{(2\pi)^3 Qi}\int_{a-i\infty}^a (2\pi Q^2 )^s \frac{\zeta(2-s)}{{\f}}\Gamma(s - \mu )\Gamma( s + \mu )\Gamma(2-s) e^{(z - \frac{3}{2}\pi i)s}ds.   \tag{7.16}
\end{align*}
\normalsize
The first is analytic for $y>-3\pi$, the second for $y>-\pi$, the third for $y>-2\pi$, the fourth is analytic on the whole complex plane by Lemma 7.1, and the next is analytic for $y<3\pi$. Therefore, (7.16) completes the proof of the continuation of $f(z,F)$ to the region $y>-\pi$.   $\square$ 

\section{Proof of Lemma 7.1}
We prove Lemma 7.1. First, we consider the case $\mu > \frac{1}{2}$. By the asymptotic expansion (5.12),
\begin{equation*}
   W_{-\frac{3}{2}, \mu}\left( \frac{e^{\frac{3}{2}\pi i - z}}{2\pi nkQ^2} \right)
   = \frac{\Gamma(2\mu)}{\Gamma(2+\mu)} \left( \frac{e^{\frac{3}{2}\pi i - z}}{2\pi nkQ^2} \right)^{\frac{1}{2}-\mu} 
   + O\left( \frac{|e^{ \left( \frac{3}{2}\pi i - z \right)\left( \frac{3}{2} -\mu \right) }|}{(2\pi nkQ^2)^{\frac{3}{2} - \mu}} \right).
\end{equation*}Hence, the inside of the curly brackets on the right hand side of (7.15) is
\begin{eqnarray*}
   \lefteqn{2\pi i(2\pi nkQ^2)^\frac{1}{2} \exp\left( \frac{e^{\frac{3}{2}\pi i - z}}{4\pi nkQ^2} +\frac{z}{2} - \frac{3}{4}\pi i \right)\Gamma(2-\mu)\Gamma(2+\mu)W_{-\frac{3}{2}, \mu}\left( \frac{e^{\frac{3}{2}\pi i - z}}{2\pi nkQ^2} \right)}
    \\
    &&-2\pi i\left\{ \sum_{m_1 =0}^M R^{(1)}_{k,n, m_1}(z,\mu) + \sum_{m_2 = 0}^{M^\prime} R^{(1)}_{k,n, m_2}(z,\mu)  \right\}
    \\
    &&= 2\pi i(2\pi nkQ^2)^\mu \exp\left( \frac{e^{\frac{3}{2}\pi i - z}}{4\pi nkQ^2} -\frac{3}{2}\mu\pi i + \mu z \right)\Gamma( 2 - \mu)\Gamma( 2\mu)
    \\
    &&+O_{Q,\mu,x} \left( \frac{1}{(2\pi nkQ^2)^{1-\mu}} \exp\left(\frac{e^{-x}}{4\pi nkQ^2} \right) \right) 
    \\
    &&- 2\pi i\left\{ \sum_{m_1 = 0}^M R^{(1)}_{k,n, m_1}(z,\mu) + \sum_{m_2 = 0}^{M^\prime} R^{(1)}_{k,n, m_2}(z,\mu)  \right\}. 
\end{eqnarray*}Now, by the Taylor expansion
\begin{align*}
   \exp\left( \frac{e^{\frac{3}{2}\pi i - z}}{4\pi nkQ^2} \right) 
   &= 1 + O\left( \left| \frac{e^{\frac{3}{2}\pi i - z}}{4\pi nkQ^2} \right| \right)
   \\
   &= 1 + O\left( \frac{e^{-x}}{4\pi nkQ^2} \right)   \tag{8.1}
\end{align*}as $n,k$ tend to infinity and (7.7), (7.8), we have
\begin{eqnarray*}
    \lefteqn{2\pi i(2\pi nkQ^2)^\frac{1}{2} \exp\left( \frac{e^{\frac{3}{2}\pi i - z}}{4\pi nkQ^2} +\frac{z}{2} - \frac{3}{4}\pi i \right)\Gamma(2-\mu)\Gamma(2+\mu)W_{-\frac{3}{2}, \mu}\left( \frac{e^{\frac{3}{2}\pi i - z}}{2\pi nkQ^2} \right)}
    \\
    &&-2\pi i\left\{ \sum_{m_1 = 0}^M R^{(1)}_{k,n, m_1}(z,\mu) + \sum_{m_2 = 0}^{M^\prime} R^{(1)}_{k,n, m_2}(z,\mu)  \right\}
    \\
    &&= O_{Q,\mu, x}\left( \frac{1}{(nk)^{1-\mu}} \right) + O_{Q,\mu,x}\left( \frac{1}{(nk)^{1 - \mu}} + \frac{1}{(nk)^{2 - \mu}} \right)
    \\
    &&-2\pi i\left\{ \sum_{m_1 = 1}^M R^{(1)}_{k,n, m_1}(z,\mu) + \sum_{m_2 = 0}^{M^\prime} R^{(1)}_{k,n, m_2}(z,\mu)  \right\}
    \\
    &&\ll O_{Q,\mu, x}\left( \frac{1}{(nk)^{1 - \mu}} \right) + \left\{ \sum_{m_1 = 1}^M \left| R^{(1)}_{k,n, m_1}(z,\mu) \right| + \sum_{m_2 = 0}^{M^\prime} \left| R^{(1)}_{k,n, m_2}(z,\mu) \right|  \right\}
    \\
    &&\ll_{Q,\mu,x} O_{Q,\mu, x}\left( \frac{1}{(nk)^{1 - \mu}} \right) + \sum_{m_1 = 1}^M \frac{1}{(nk)^{m_1 - \mu}} + \sum_{m_2 = 0}^{M^\prime} \frac{1}{(nk)^{m_2 + \mu}}
    \\
    &&= O_{Q,\mu,x}\left( \frac{1}{(nk)^{1 - \mu}} \right)
\end{eqnarray*}Hence, $I_1 (z,F)$ is evaluated as follows :
\begin{align*}
   I_1 (z,F)
   &\ll_{Q,\mu,x} \sum_{k,n=1}^\infty \frac{|{\mf}|}{kn^2} \cdot O_{Q,\mu, x}\left( \frac{1}{(nk)^{1 - \mu}} \right)   
   \\
   &= O_{Q,\mu,x}\left(\sum_{k,n=1}^\infty \frac{|{\mf}|}{k^{2 - \mu}n^{3 - \mu}} \right).
\end{align*}Therefore, the series on the right hand side in (7.15) is convergent for $\frac{1}{2} < \mu < 1$. 

Secondly, in the case $\mu = \frac{1}{2}$, by (5.13),
\begin{align*}
  W_{-\frac{3}{2},\frac{1}{2}}\left( \frac{e^{\frac{3}{2}\pi i - z}}{2\pi nkQ^2} \right)
   &= \frac{1}{\Gamma\left( \frac{5}{2} \right)} + O_y \left( \frac{e^{-x}}{2\pi nkQ^2}\log \left( \frac{e^{-x}}{2\pi nkQ^2} \right) \right)
   \\
   &= \frac{1}{\Gamma\left( \frac{5}{2} \right)} + O_y \left( \left( \frac{e^{-x}}{2\pi nkQ^2} \right)^{1-\delta} \right),
\end{align*}where $\delta$ is any positive real number. By the same calculation as in the first case and using (7.10), the inside of the curly brackets on the right hand side of (7.15) is evaluated as follows :
\begin{eqnarray*}
   \lefteqn{2\pi i(2\pi nkQ^2)^\frac{1}{2} \exp\left( \frac{e^{\frac{3}{2}\pi i - z}}{4\pi nkQ^2} +\frac{z}{2} - \frac{3}{4}\pi i \right)\Gamma\left( \frac{3}{2} \right)\Gamma\left( \frac{5}{2} \right)W_{-\frac{3}{2}, \frac{1}{2}}\left( \frac{e^{\frac{3}{2}\pi i - z}}{2\pi nkQ^2} \right)}
   \\
   &&- 2\pi i\sum_{m= 0}^{M^{\prime\prime}} R^{(1)}_{k,n, m,\frac{1}{2}}(z)
   \\
    &&\ll O_{Q, x,y}\left( \frac{1}{(nk)^{\frac{1}{2}-\delta}} \right) + \sum_{m=1}^{M^{\prime\prime}} \frac{\log nk}{(nk)^{m - \frac{1}{2}}} 
    \\
    &&= O_{Q,x,y,M^{\prime\prime}}\left(\frac{1}{(nk)^{\frac{1}{2}-\delta}}\right).
\end{eqnarray*}Hence, $I_1 (z,F)$ is evaluated as follows :
\begin{align*}
I_1 (z,F)
&\ll \sum_{k,n=1}^\infty \frac{|\mu_F (k)|}{kn^2} \cdot O_{Q,x,y,M^{\prime\prime}}\left( \frac{1}{(nk)^{\frac{1}{2}-\delta}}\right)
\\
&= O_{Q,x,y,M^{\prime\prime}} \left( \sum_{k,n=1}^\infty \frac{|\mu_F (k)|}{k^{\frac{3}{2}-\delta}n^{\frac{5}{2}-\delta}} \right).
\end{align*}
Therefore, the series on the right hand side in (7.15) is convergent for $\mu = \frac{1}{2}$. 

Thirdly, in the case $0 < \mu < \frac{1}{2}$, by (5.14),
\begin{equation*}
   W_{-\frac{3}{2}, \mu}\left( \frac{e^{\frac{3}{2}\pi i - z}}{2\pi nkQ^2} \right)
   = \frac{\Gamma(2\mu)}{\Gamma(2+\mu)} \left( \frac{e^{\frac{3}{2}\pi i - z}}{2\pi nkQ^2} \right)^{\frac{1}{2}-\mu} 
   + O\left( \frac{e^{-x\left( \mu + \frac{1}{2} \right)}}{(2\pi nkQ^2)^{ \mu + \frac{1}{2}}} \right)
\end{equation*}and
\begin{eqnarray*}
   \lefteqn{2\pi i(2\pi nkQ^2)^\frac{1}{2} \exp\left( \frac{e^{\frac{3}{2}\pi i - z}}{4\pi nkQ^2} +\frac{z}{2} - \frac{3}{4}\pi i \right)\Gamma(2-\mu)\Gamma(2+\mu)W_{-\frac{3}{2}, \mu}\left( \frac{e^{\frac{3}{2}\pi i - z}}{2\pi nkQ^2} \right)}
    \\
    &&-2\pi i\left\{ \sum_{m_1 = 0}^M R^{(1)}_{k,n, m_1}(z,\mu) + \sum_{m_2 = 0}^{M^\prime} R^{(1)}_{k,n, m_2}(z,\mu)  \right\}
    \\
    &&\ll_{Q,\mu,x} O_{Q,\mu, x}\left( \frac{1}{(nk)^{1 - \mu}} \right) + \sum_{m_1 = 1}^M \frac{1}{(nk)^{m_1 - \mu}} + \sum_{m_2 = 0}^{M^\prime} \frac{1}{(nk)^{m_2 + \mu}}
    \\
    &&= O_{Q,\mu, x}\left( \frac{1}{(nk)^{1 - \mu}} \right).
\end{eqnarray*}Therefore, the series on the right hand side in (7.15) is convergent for $0 < \mu < \frac{1}{2}$.

Finally, in the case $\mu=0$, by (5.16),
\begin{align*}
   W_{-\frac{3}{2},0}\left( \frac{e^{\frac{3}{2}\pi i - z}}{2\pi nkQ^2} \right)
   &= -\frac{e^{\frac{3}{4}\pi i -\frac{z}{2}}}{(2\pi nkQ^2)^\frac{1}{2}}\left\{ \log\left( \frac{e^{-x}}{2\pi nkQ^2} \right) + \left( \frac{3}{2}\pi - y \right)i + \frac{\Gamma^\prime}{\Gamma}(2) + 2C_0 \right\}
   \\
   &+ O_{Q,x} \left( \frac{1}{(nk)^\frac{3}{2}}\log\left( \frac{e^{-x}}{2\pi nkQ^2} \right) \right).
\end{align*}Using the recurrence formula
\begin{equation*}
   \psi(s+1) = \frac{1}{s} + \psi(s)   \tag{8.2}
\end{equation*}(see~\cite{L}) and (7.9), the inside of the curly brackets on the right hand side of (7.15) is evaluated as follows :
\begin{eqnarray*}
   \lefteqn{2\pi i (2\pi nkQ^2)^\frac{1}{2} \exp\left( \frac{e^{\frac{3}{2}\pi i - z} }{4\pi nkQ^2} + \frac{z}{2} - \frac{3}{4}\pi i \right) \Gamma(2)^2 W_{ -\frac{3}{2},0 } \left( \frac{e^{\frac{3}{2}\pi i - z} }{2\pi nkQ^2} \right) - 2\pi i\sum_{m = 0}^{ M^{\prime\prime} } R^{(1)}_{k,n, m,0}(z)}
   \\
   &&= -2\pi i \left( \frac{\Gamma^\prime}{\Gamma}(2) + C_0 -1 \right) + O_{Q,x,y} \left( \frac{1}{nk}\log \left( \frac{e^{-x}}{2\pi nkQ^2} \right) \right) -2\pi i\sum_{m = 1}^{ M^{\prime\prime} } R^{(1)}_{k,n,m,0}(z)
    \\
    &&= O_{Q,x,y} \left( \frac{1}{nk}\log \left( \frac{e^{-x}}{2\pi nkQ^2} \right) \right) - 2\pi i\sum_{m = 1}^{ M^{\prime\prime} } R^{(1)}_{k,n,m,0}(z)
    \\
    &&= O_{Q,x,y} \left( \frac{1}{nk} \cdot \frac{1}{(nk)^{1-\delta}} \right) - 2\pi i\sum_{m = 1}^{ M^{\prime\prime} } R^{(1)}_{k,n,m,0}(z)
    \\
    &&= O_{Q,x,y}\left( \frac{1}{(nk)^{2-\delta}} \right) -2\pi i\sum_{m = 1}^{M^{\prime\prime}} R^{(1)}_{k,n,m,0}(z)
    \\
    &&\ll_{Q,x,y}  O_{Q,x,y} \left( \frac{1}{(nk)^{2-\delta}} \right) + \sum_{m = 1}^{M^{\prime\prime}} \frac{\log nk}{(nk)^m}
    \\
    &&\ll O_{M^{\prime\prime}}\left( \frac{1}{(nk)^{1-\delta}} \right),
\end{eqnarray*}where $\delta$ is any positive real number. Hence, $I_1 (z,F)$ is evaluated as follows :
\begin{align*}
   I_1 (z,F)
   &\ll_{M^{\prime\prime}}  \sum_{k,n=1}^\infty \frac{|{\mf}|}{kn^2} \cdot \frac{1}{(nk)^{1-\delta}}
   \\
   &= \sum_{k,n=1}^\infty \frac{|{\mf}|}{k^{2-\delta}n^{3-\delta}}.
\end{align*}Therefore, the series on the right hand side in (7.15) is convergent for $\mu = 0$.  In summary, 
\begin{align*}
   I_1 (z,F) \leq \sum_{k,n=1}^\infty \frac{|{\mf}|}{kn^2} \cdot \max\left\{ \frac{1}{(nk)^{ 1 - \mu}},\frac{1}{(nk)^{\frac{1}{2}-\delta}} \right\}.  \tag{8.3}
\end{align*}By (5.2),(5.3) and (5.4), the Whittaker function $W_{-\frac{3}{2},\mu} \left( \frac{e^{\frac{3}{2}\pi i-z} }{2\pi nkQ^2} \right)$ is analytic for all $z \in \mathbb{C}$. Therefore, We have the desired result. $\square$ 

\section{Proof of Corollary 4.3 and another proof}
We prove Corollary 4.3. We use the following lemma similar to Lemma 3.1. We can prove this lemma by modifying the proof of Lemma 3.1.
\begin{lemma}
Let $F \in \mathcal{S}$ and let $T$ be sufficiently large. Moreover, let $H=D\log\log T$ be fixed, where $D$ is a large positive constant. In any subinterval of length $1$ in $[-T-H, -T+H]$ there are lines $t=t_0$ such that
\begin{equation*}
   |F( \sigma+it_0 )|^{-1} = O(\exp(C(\log T)^2))   \tag{9.1}  
\end{equation*}uniformly in $\sigma \geq -2$.
\end{lemma}
We consider the integral
\begin{equation*}
   \int_{\mathscr{L}^\prime} \frac{\zeta(s-1)}{F(s)}e^{zs}ds,  \tag{9.2}
\end{equation*}where $\mathscr{L}^\prime$ is the contour symmetrical upon the real axis to $\mathscr{L}$ in (3.3). By Lemma 9.1, the integral along the lower side of the contour tends to $0$ as $n$ tends to infinity for $z \in \mathbb{H}^{-}$. Then, we have by residue theorem and the definition (4.9), in a similar manner as (4.1),
\begin{equation*}
   2\pi i f^{-}(z,F) = f_{1}^{-}(z,F) +  f_{2}^{-}(z,F) +  f_{3}^{-}(z,F),   \tag{9.3}
\end{equation*}where
\begin{equation*}
   f_{1}^{-}(z,F) = \int_{a}^{a-i\infty} \frac{\zeta(s-1)}{F(s)}e^{sz}ds   \tag{9.4} 
\end{equation*}is analytic on $\mathbb{H}^-$, 
\begin{equation*}
  f_{2}^{-}(z,F) = \int_{\overline{L}} \frac{\zeta(s-1)}{F(s)}e^{sz}ds   \tag{9.5}
\end{equation*}is analytic on the whole complex plane. We consider the same setting as in $L$ for the curve $\overline{L}$. In the same way as obtaining (6.1),
\begin{align*}
   f_{3}^{-}(z,F)
   &= \int_{b - i\infty}^b \left( \sum_{n=1}^\infty \frac{g(n)}{n^s} \right)e^{zs}ds
   \\
   &= \sum_{n=1}^\infty g(n) \int_{b - i\infty}^b e^{s(z-\log n)}ds
   \\
   &= e^{bz}\sum_{n=1}^\infty \frac{g(n)}{n^b (z-\log n)}   \tag{9.6}
\end{align*}is meromorphic on the whole complex plane. Now we already know that $f_{1}^{-} (z,F)$ is analytic for $y<0$, and we have to continue to $y<\pi$ just as in the case of $f_1 (z,F)$ (see Section 7). By the functional equation for $\zeta(s)$ and $F(s)$, we have
\begin{align*}
  f_{1}^{-} (z,F)
   &= \int_{a}^{a-i\infty} \frac{\zeta(s-1)}{F(s)}e^{zs}ds
   \\
   &= -\int_{a-i\infty}^a \frac{\zeta(s-1)}{F(s)}e^{zs}ds
   \\
   &= f_{11}^- (z,F) + f_{12}^- (z,F) + f_{13}^- (z,F) + f_{14}^- (z,F),  \tag{9.7}
\end{align*}where
\begin{equation*}
   f_{11}^{-}(z,F) 
   = \frac{\overline{\omega} e^{-\mu\pi i} }{(2\pi)^3 Qi}\int_{a-i\infty}^a (2\pi Q^2 )^s \frac{\zeta(2-s)}{{\f}}\Gamma(s+\mu)\Gamma(s-\mu)\Gamma(2-s)e^{\left( z + \frac{3}{2}\pi i \right)s}ds   \tag{9.8}
\end{equation*}is analytic for $y < 0$,
\begin{equation*}
   f_{12}^{-}(z,F)
   = -\frac{\overline{\omega} e^{\mu\pi i} }{(2\pi)^3 Qi}\int_{a-i\infty}^a (2\pi Q^2 )^s \frac{\zeta(2-s)}{{\f}}\Gamma(s+\mu)\Gamma(s-\mu)\Gamma(2-s)e^{\left( z-\frac{\pi}{2}i \right)s}ds    \tag{9.9}
\end{equation*}for $y < 2\pi$, 
\begin{equation*}
   f_{13}^{-}(z,F)
  =  \frac{\overline{\omega} e^{-\mu\pi i} }{(2\pi)^3 Qi} \int_{a-i\infty}^a (2\pi Q^2 )^s \frac{\zeta(2-s)}{{\f}}\Gamma(s+\mu)\Gamma(s-\mu)\Gamma(2-s)e^{\left( z+\frac{\pi}{2}i \right)s}ds    \tag{9.10}
\end{equation*}for $y < \pi$,
\begin{equation*}
  f_{14}^{-}(z,F)
  = -\frac{\overline{\omega} e^{\mu\pi i}}{(2\pi)^3 Qi} \int_{a-i\infty}^a (2\pi Q^2 )^s \frac{\zeta(2-s)}{{\f}}\Gamma(s+\mu)\Gamma(s-\mu)\Gamma(2-s)e^{\left( z-\frac{3}{2}\pi i \right)s}ds   \tag{9.11}
\end{equation*} for $y < 3\pi$. Splitting the integral on the right hand side in (8.8) just as in the case of $f_{14}(z,F)$, we have
\[ f_{11}^{-}(z,F) = I_{1}^{-} (z,F) + I_{2}^{-} (z,F), \]where
\begin{equation*}
    I_{1}^{-} (z,F) = \frac{\overline{\omega}e^{-\mu\pi i}}{(2\pi)^3 Qi} \int_{a-i\infty}^{a+i\infty}(2\pi Q^2)^s \frac{\zeta(2-s)}{{\f}}\Gamma(s+\mu)\Gamma(s-\mu)\Gamma(2-s)e^{\left( z + \frac{3}{2}\pi i \right)s}ds   \tag{9.12}
\end{equation*}and
\begin{equation*}
    I_{2}^{-} (z,F) = -\frac{\overline{\omega}e^{-\mu\pi i}}{(2\pi)^3 Qi} \int_{a}^{a+i\infty}(2\pi Q^2)^s \frac{\zeta(2-s)}{{\f}}\Gamma(s+\mu)\Gamma(s-\mu)\Gamma(2-s)e^{\left( z + \frac{3}{2}\pi i \right)s}ds.   \tag{9.13}
 \end{equation*}We see that the integral $I_{2}^{-} (z,F)$ is convergent for $y > -3\pi$ by the same way as in (7.5). Since $f_{11}^{-}(z,F)$ is analytic for $y < 0$, the integral $I_{1}^{-} (z,F)$ is convergent for $-3\pi < y < 0$ and we can calculate $I_{1}^{-}(z,F)$ for $-3\pi < y < 0$ in a similar way as (7.4) (Section 7). Let $m_1, m_2$ and $m$ be non-negative integers. By taking the path of integration $C$ in (7.12), we have for $\mu$ which is not integer except zero and $|y+\frac{3}{2}\pi| < \frac{3}{2}\pi$
\footnotesize
\begin{align*}
    I_{1}^{-}(z,F) 
    &= \frac{\overline{\omega}e^{-\mu\pi i}}{(2\pi)^3 Qi} \sum_{k,n=1}^\infty \frac{{\mf}}{kn^2} \times 2\pi i \left\{(2\pi nkQ^2)^\frac{1}{2}\exp\left( \frac{3}{4}\pi i + \frac{z}{2} + \frac{e^{-\frac{3}{2}\pi i-z} }{4\pi nkQ^2} \right) \right.
    \\
    &\times \left. \Gamma(2+\mu)\Gamma(2-\mu)W_{-\frac{3}{2}, \> \mu}\left( \frac{e^{-\frac{3}{2}\pi i - z}}{2\pi nkQ^2} \right) - \sum_{m_1 = 0}^M R^{(2)}_{k,n, m_1}(z,\mu) \hspace{-0.01cm} - \hspace{-0.01cm} \sum_{m_2 = 0}^{M^\prime} R^{(2)}_{k,n, m_2}(z,\mu) \hspace{-0.01cm} \right\},  \tag{9.14}
\end{align*}
\normalsize 
where
\begin{gather*}
   R^{(2)}_{k,n, m_1}(z,\mu)
    = \frac{(-1)^{m_1}}{m_{1}!}\Gamma(2\mu-m_1)\Gamma(2-\mu+m_1)(2\pi nkQ^2)^{\mu-m_1}e^{\left( z + \frac{3}{2}\pi i \right)( \mu - m_1 )},  \tag{9.15}
    \\
    R^{(2)}_{k,n, m_2}(z,\mu)
    = \frac{(-1)^{m_2}}{m_{2}!}\Gamma(-2\mu-m_2)\Gamma(2+\mu+m_2)(2\pi nkQ^2)^{-\mu-m_2}e^{\left( z + \frac{3}{2}\pi i \right)( -\mu - m_2 )},  \tag{9.16}
    \\
    R^{(2)}_{k,n,m,0}(z)
    =  \frac{m+1}{m!}\frac{e^{-\left( z+\frac{3}{2}\pi i \right)m}}{(2\pi nkQ^2)^m}\left\{ \log(2\pi nkQ^2) + z  + \frac{3}{2}\pi i  + \sum_{k_1 = 1}^m \frac{1}{k_1} - C_0 - \frac{1}{m+1} \right\}  \tag{9.17}
\end{gather*}and
\begin{align*}
   R^{(2)}_{k,n,m,\frac{1}{2}}(z)
   &= \frac{\Gamma\left( \frac{3}{2}+m \right)}{(m!)^2}\frac{e^{\left( \frac{1}{2}-m \right)\left( z + \frac{3}{2}\pi i \right)}}{(2\pi nkQ^2)^{m-\frac{1}{2}}}\left\{ m\left( \psi\left( \frac{3}{2}+m \right) -2\psi(m+1) \right)  \right.
   \\
   &- \left. m\left( \log (2\pi nkQ^2) + z + \frac{3}{2}\pi i \right) + 1 \right\}      \tag{9.18}
\end{align*}
 are residues of the integrand in (9.12) at  $s = \mu - m_1, -\mu - m_2, -m$ and $\frac{1}{2}-m$ respectively. The convergence of the series on the right hand side in (9.14) follows in a similar manner as the consideration in (7.15). Finally, by (9.7)-(9.18) we obtain the following continuation of $f_{1}^{-}(z,F)$ to $y<\pi$ : For $F \in \mathcal{S}^{\text{poly}}$ whose $(r,\lambda_j) = (1,1)$ for all $j$ in (1.8) and $0 \leq \mu < 1$, 
\footnotesize
\begin{align*}
    f_{1}^{-} (z,F)
    &= \frac{\overline{\omega}e^{-\mu\pi i}}{(2\pi)^3 Qi} \sum_{k,n=1}^\infty \frac{{\mf}}{kn^2}  \times 2\pi i \left\{ (2\pi nkQ^2)^\frac{1}{2}\exp\left( \frac{3}{4}\pi i + \frac{z}{2} + \frac{e^{-\frac{3}{2}\pi i-z} }{4\pi nkQ^2} \right)  \right.
    \\
    &\times \left. \Gamma(2+\mu)\Gamma(2-\mu)W_{-\frac{3}{2}, \> \mu}\left( \frac{e^{-\frac{3}{2}\pi i - z}}{2\pi nkQ^2} \right) - \sum_{m_1 = 0}^M R^{(2)}_{k,n, m_1}(z,\mu) - \sum_{m_2 = 0}^{M^\prime} R^{(2)}_{k,n, m_2}(z,\mu) \right\}
    \\
    &- \frac{\overline{\omega}e^{-\overline{\mu}\pi i}}{(2\pi)^3 Qi}\int_{a}^{a+i\infty} (2\pi Q^2)^s \frac{\zeta(2-s)}{{\f}}\Gamma(s+\mu)\Gamma(s-\mu)\Gamma(2-s)e^{\left( z + \frac{3}{2}\pi i \right)s}ds
    \\
    &-\frac{\overline{\omega} e^{\mu\pi i} }{(2\pi)^3 Qi}\int_{a-i\infty}^a (2\pi Q^2 )^s \frac{\zeta(2-s)}{{\f}}\Gamma(s+\mu)\Gamma(s-\mu)\Gamma(2-s)e^{\left( z-\frac{\pi}{2}i \right)s}ds
    \\
    &+ \frac{\overline{\omega} e^{-\mu\pi i} }{(2\pi)^3 Qi} \int_{a-i\infty}^a (2\pi Q^2 )^s \frac{\zeta(2-s)}{{\f}}\Gamma(s+\mu)\Gamma(s-\mu)\Gamma(2-s)e^{\left( z+\frac{\pi}{2}i \right)s}ds
    \\
    &-\frac{\overline{\omega} e^{\mu\pi i}}{(2\pi)^3 Qi} \int_{a-i\infty}^a (2\pi Q^2 )^s \frac{\zeta(2-s)}{{\f}}\Gamma(s+\mu)\Gamma(s-\mu)\Gamma(2-s)e^{\left( z-\frac{3}{2}\pi i \right)s}ds.  \tag{9.19}
\end{align*}
\normalsize
Since a lemma similar to Lemma7.1 holds for the series on the right hand side in (9.19), the series we now consider is also absolutely and uniformly convergent on every compact subset on the whole complex plane. Therefore, we complete the continuation of $f^{-}(z,F)$ analytic for $y<0$ to the region $y<\pi$.

Also, Corollary4.3 can be proved form Theorem4.2 and the definition (4.9) directly as follows : For $z \in \mathbb{H}^-$ and $\rho$\hspace{0.01cm} with $\text{Im} \hspace{0.1cm} \rho <0$, 
\begin{align*}
   f^- (z,F)
   &= \sum_{\rho} \frac{e^{\rho z}\zeta(\rho-1)}{F^\prime (\rho)} 
   \\
   &=  \overline{\sum_{\rho} \frac{ \overline{\zeta(\rho-1)} }{\overline{F^\prime (\rho)} }e^{\overline{\rho z} } }
   \\
   &=  \overline{\sum_{\rho^\prime } \frac{\zeta(\rho^\prime -1)}{\overline{F^\prime \left(\overline{\rho^\prime } \right) } }e^{\rho^\prime \overline{z} } },
\end{align*}where $\rho^\prime = \overline{\rho}$. We recall the definition $\overline{F} (s) = \overline{F(\overline{s})}$ for $F \in {\se}$. Hence, the sum on the right hand side in the third line yields
\begin{equation*}
   \overline{\sum_{\rho^\prime} \frac{\zeta(\rho^\prime -1)}{\overline{F^\prime}(\rho^\prime)  }e^{\rho^\prime \overline{z} } }
   = f(\overline{z},\overline{F}).
\end{equation*}Here, we use the fact that if $F \in {\s}$, then so is $\overline{F} \in {\s}$. Of course, ${\s}$ may be replaced by ${\se}$. By Theorem 4.2, the function $f(z,F)$ has a meromorphic continuation to $y > -\pi$. Hence $f(\overline{z},\overline{F})$ has a meromorphic continuation to $y < \pi$.

\section{Proof of Theorem 4.4}
We add (7.16) to (9.19). Since some integrals are canceled, we have for $|y| < \pi$
\footnotesize
\begin{align*}
   f_1 (z,F) + f_{1}^{-}(z,F)
   &= f_{11}(z,F) + f_{12}(z,F) + f_{13}(z,F) + \frac{\overline{\omega}e^{\mu\pi i}}{(2\pi)^3 Qi}I_1 (z,F) +  \frac{\overline{\omega}e^{\mu\pi i}}{(2\pi)^3 Qi}I_2 (z,F) 
   \\
   &+ I_{1}^- (z,F) + I_{2}^- (z,F) + f_{12}^- (z,F) + f_{13}^- (z,F) + f_{14}^- (z,F)  
   \\ 
   &= \frac{\overline{\omega}e^{\mu\pi i}i}{(2\pi)^3 Q} \sum_{k,n=1}^\infty \frac{{\mf}}{kn^2}  \times 2\pi i \left\{ (2\pi nkQ^2)^\frac{1}{2}\exp\left( -\frac{3}{4}\pi i + \frac{z}{2} + \frac{e^{\frac{3}{2}\pi i-z} }{4\pi nkQ^2} \right)  \right.
   \\
   &\times \left. \Gamma(2+\mu)\Gamma(2-\mu)W_{-\frac{3}{2}, \> \mu}\left( \frac{e^{\frac{3}{2}\pi i - z}}{2\pi nkQ^2} \right) - \sum_{m_1 = 0}^M R^{(1)}_{k,n, m_1}(z,\mu) - \sum_{m_2 = 0}^{M^\prime} R^{(1)}_{k,n, m_2}(z,\mu) \right\}
  \\
  &+\frac{\overline{\omega}e^{-\mu\pi i}}{(2\pi)^3 Qi} \sum_{k,n=1}^\infty \frac{{\mf}}{kn^2} \times 2\pi i \left\{ (2\pi nkQ^2)^\frac{1}{2} \exp\left( \frac{3}{4}\pi i + \frac{z}{2} + \frac{e^{-\frac{3}{2}\pi i-z} }{4\pi nkQ^2} \right) \right.
  \\
  &\times \left. \Gamma(2+\mu)\Gamma(2-\mu) W_{-\frac{3}{2}, \> \mu}\left( \frac{e^{-\frac{3}{2}\pi i - z}}{2\pi nkQ^2} \right) - \sum_{m_1 = 0}^M R^{(2)}_{k,n, m_1}(z,\mu) - \sum_{m_2 = 0}^{M^\prime} R^{(2)}_{k,n, m_2}(z,\mu) \right\}
   \\
   &+ A_{1}(z,F) + A_{2}(z,F),
\end{align*}
\normalsize
where
\begin{equation*}
   A_{1} (z,F) = -\frac{\overline{\omega} e^{\mu\pi i}}{(2\pi)^3 Qi} \int_{a-i\infty}^{a+i\infty} (2\pi Q^2 )^s \frac{\zeta(2-s)}{{\f}}\Gamma(s+\mu)\Gamma(s-\mu)\Gamma(2-s)e^{\left( z - \frac{\pi}{2}i \right)s}ds   \tag{10.1}   
\end{equation*}and
\begin{equation*}
   A_{2} (z,F) = \frac{\overline{\omega} e^{-\mu\pi i}}{(2\pi)^3 Qi} \int_{a-i\infty}^{a+i\infty} (2\pi Q^2 )^s \frac{\zeta(2-s)}{{\f}}\Gamma(s+\mu)\Gamma(s-\mu)\Gamma(2-s)e^{\left( z + \frac{\pi}{2}i \right)s}ds.  \tag{10.2}
\end{equation*}The integrals $A_1 (z,F)$ and $A_2 (z,F)$ are convergent for $|y| < \pi$ and we can obtain series expressions of them involving Whittaker functions in a way similar to the case of $I_1 (z,F)$ and $I_{1}^{-}(z,F)$. We have for $|y|<\pi, F \in \mathcal{S}^{\text{poly}}$ whose $(r,\lambda_j) = (1,1)$ for all $j$ in (1.8) and $0 \leq \mu < 1$
\begin{align*}
    A_1 (z,F)
    &= -\frac{\overline{\omega}e^{\mu\pi i}}{(2\pi)^3 Qi} \sum_{k,n=1}^\infty \frac{{\mf}}{kn^2} \times 2\pi i
    \left\{ (2\pi nkQ^2)^\frac{1}{2}\exp\left( -\frac{\pi}{4}i + \frac{z}{2} + \frac{e^{\frac{\pi}{2}i-z} }{4\pi nkQ^2} \right)  \right.
    \\
    &\times \left. \Gamma(2+\mu)\Gamma(2-\mu)W_{-\frac{3}{2}, \> \mu} \left( \frac{e^{\frac{\pi}{2}i - z}}{2\pi nkQ^2} \right) -\sum_{m_1 = 0}^M R^{(3)}_{k,n, m_1}(z,\mu) - \sum_{m_2 = 0}^{M^\prime} R^{(3)}_{k,n, m_2}(z,\mu) \right\},
\end{align*}
\normalsize
where
\begin{gather*}
   R^{(3)}_{k,n, m_1}(z,\mu)
    = \frac{(-1)^{m_1}}{m_{1}!}\Gamma(2\mu-m_1)\Gamma(2-\mu+m_1)(2\pi nkQ^2)^{\mu-m_1}e^{\left( z - \frac{\pi}{2}i \right)( \mu - m_1 )},  \tag{10.3}
    \\
    R^{(3)}_{k,n, m_2}(z,\mu)
    = \frac{(-1)^{m_2}}{m_{2}!}\Gamma(-2\mu-m_2)\Gamma(2+\mu+m_2)(2\pi nkQ^2)^{-\mu-m_2}e^{\left( z - \frac{\pi}{2}i \right)( -\mu - m_2 )},  \tag{10.4}
    \\
    R^{(3)}_{k,n,m,0}(z)
    =  \frac{m+1}{m!}\frac{e^{-\left( z-\frac{\pi}{2}i \right)m}}{(2\pi nkQ^2)^m}\left\{ \log(2\pi nkQ^2) + z - \frac{\pi}{2}i + \sum_{k_1 = 1}^m \frac{1}{k_1} - C_0 - \frac{1}{m+1} \right\}   \tag{10.5}
\end{gather*}and
\begin{align*}
   R^{(3)}_{k,n,m,\frac{1}{2}}(z)
   &= \frac{\Gamma\left( \frac{3}{2}+m \right)}{(m!)^2}\frac{e^{\left( \frac{1}{2}-m \right)\left( z -\frac{\pi}{2} i \right)}}{(2\pi nkQ^2)^{m-\frac{1}{2}}}\left\{ m\left( \psi\left( \frac{3}{2}+m \right) -2\psi(m+1) \right)  \right.
   \\
   &- \left. m\left( \log (2\pi nkQ^2) + z - \frac{\pi}{2}i \right) + 1 \right\} \tag{10.6}
\end{align*}are residues of the integrand in $A_1(z,F)$ at  $s = \mu - m_1, -\mu - m_2, -m$ and $\frac{1}{2}-m$ respectively. Similarly, for $|y|<\pi, F \in \mathcal{S}^{\text{poly}}$ whose $(r,\lambda_j) = (1,1)$ for all $j$ in (1.8) and $0 \leq \mu < 1$
\begin{align*}
    A_2 (z,F) 
    &= \frac{\overline{\omega}e^{-\mu\pi i}}{(2\pi)^3 Qi} \sum_{k,n=1}^\infty \frac{{\mf}}{kn^2}  \times 2\pi i \left\{ (2\pi nkQ^2)^\frac{1}{2}\exp\left( \frac{\pi}{4}i + \frac{z}{2} + \frac{e^{-\frac{\pi}{2}i-z} }{4\pi nkQ^2} \right)  \right.
    \\
    &\times \left. \Gamma(2+\mu)\Gamma(2-\mu)W_{-\frac{3}{2}, \> \mu}\left( \frac{e^{-\frac{\pi}{2}i - z}}{2\pi nkQ^2} \right) - \sum_{m_1 = 0}^M R^{(4)}_{k,n, m_1}(z,\mu) - \sum_{m_2 = 0}^{M^\prime} R^{(4)}_{k,n, m_2}(z,\mu) \right\},
\end{align*}
\normalsize
where 
\begin{gather*}
    R^{(4)}_{k,n, m_1}(z,\mu)
    = \frac{(-1)^{m_1}}{m_{1}!}\Gamma(2\mu-m_1)\Gamma(2-\mu+m_1)(2\pi nkQ^2)^{\mu-m_1}e^{\left( z + \frac{\pi}{2}i \right)( \mu - m_1 )},  \tag{10.7}
    \\
    R^{(4)}_{k,n, m_2}(z,\mu)
    = \frac{(-1)^{m_2}}{m_{2}!}\Gamma(-2\mu-m_2)\Gamma(2+\mu+m_2)(2\pi nkQ^2)^{-\mu-m_2}e^{\left( z + \frac{\pi}{2}i \right)( -\mu - m_2 )},  \tag{10.8}
    \\
    R^{(4)}_{k,n,m,0}(z)
    =  \frac{m+1}{m!}\frac{e^{-\left( z+\frac{\pi}{2}i \right)m}}{(2\pi nkQ^2)^m}\left\{ \log(2\pi nkQ^2) + z  + \frac{\pi}{2}i + \sum_{k_1 = 1}^m \frac{1}{k_1} - C_0 - \frac{1}{m+1} \right\}  \tag{10.9}    
\end{gather*} and
\begin{align*}
   R^{(4)}_{k,n,m,\frac{1}{2}}(z) 
   &= \frac{\Gamma\left( \frac{3}{2}+m \right)}{(m!)^2}\frac{e^{\left( \frac{1}{2}-m \right)\left( z + \frac{\pi}{2} i \right)}}{(2\pi nkQ^2)^{m-\frac{1}{2}}}\left\{ m\left( \psi\left( \frac{3}{2}+m \right) -2\psi(m+1) \right)  \right.
   \\
   &- \left. m\left( \log (2\pi nkQ^2) + z + \frac{\pi}{2}i \right) + 1 \right\}   
     \tag{10.10}  
\end{align*}are residues of the integrand in $A_2 (z,F)$ at  $s = \mu - m_1, -\mu - m_2, -m$ and $\frac{1}{2}-m$ respectively. Finally, for $|y|<\pi, F \in \mathcal{S}^{\text{poly}}$ whose $(r,\lambda_j) = (1,1)$ for all $j$ in (1.8) and $0 \leq \mu < 1$, we have the series expression for $ f_1 (z,F) + f_{1}^{-}(z,F)$
\footnotesize
\begin{align*}
   f_1 (z,F) + f_{1}^{-}(z,F)
   &= \frac{\overline{\omega}e^{\mu\pi i}i}{(2\pi)^3 Q} \sum_{k,n=1}^\infty \frac{{\mf}}{kn^2} \times 2\pi i \left\{ (2\pi nkQ^2)^\frac{1}{2}\exp\left( -\frac{3}{4}\pi i + \frac{z}{2} + \frac{e^{\frac{3}{2}\pi i - z}}{4\pi nkQ^2} \right) \right. 
   \\
   &\times \left. \Gamma(2+\mu)\Gamma(2-\mu)W_{-\frac{3}{2}, \> \mu}\left( \frac{e^{\frac{3}{2}\pi i - z}}{2\pi nkQ^2} \right) - \sum_{m_1 = 0}^M R^{(1)}_{k,n, m_1}(z,\mu) - \sum_{m_2 = 0}^{M^\prime} R^{(1)}_{k,n, m_2}(z,\mu) \right\}
   \\
   &+\frac{\overline{\omega}e^{-\mu\pi i}}{(2\pi)^3 Qi} \sum_{k,n=1}^\infty \frac{{\mf}}{kn^2} \times 2\pi i \left\{ (2\pi nkQ^2)^\frac{1}{2}\exp\left( \frac{3}{4}\pi i + \frac{z}{2} + \frac{e^{-\frac{3}{2}\pi i - z} }{4\pi nkQ^2} \right)  \right.
   \\
    &\times \left. \Gamma(2+\mu)\Gamma(2-\mu)W_{-\frac{3}{2}, \> \mu}\left( \frac{e^{-\frac{3}{2}\pi i - z}}{2\pi nkQ^2} \right) - \sum_{m_1 = 0}^M R^{(2)}_{k,n, m_1}(z,\mu) - \sum_{m_2 = 0}^{M^\prime} R^{(2)}_{k,n, m_2}(z,\mu) \right\} 
    \\
    &-\frac{\overline{\omega}e^{\mu\pi i}}{(2\pi)^3 Qi} \sum_{k,n=1}^\infty \frac{{\mf}}{kn^2}  \times 2\pi i \left\{ (2\pi nkQ^2)^\frac{1}{2}\exp\left( -\frac{\pi}{4}i + \frac{z}{2} + \frac{e^{\frac{\pi}{2}i - z} }{4\pi nkQ^2} \right) \right. 
    \\
    &\times \left. \Gamma(2+\mu)\Gamma(2-\mu)W_{-\frac{3}{2}, \> \mu}\left( \frac{e^{\frac{\pi}{2}i - z}}{2\pi nkQ^2} \right) - \sum_{m_1 = 0}^M R^{(3)}_{k,n, m_1}(z,\mu) - \sum_{m_2 = 0}^{M^\prime} R^{(3)}_{k,n, m_2}(z,\mu) \right\}
    \\
    &+\frac{\overline{\omega}e^{-\mu\pi i}}{(2\pi)^3 Qi} \sum_{k,n=1}^\infty \frac{{\mf}}{kn^2}  \times 2\pi i \left\{ (2\pi nkQ^2)^\frac{1}{2} \exp\left( \frac{\pi}{4}i + \frac{z}{2} + \frac{e^{-\frac{\pi}{2}i - z} }{4\pi nkQ^2} \right) \right. 
    \\
    &\times \left. \Gamma(2+\mu)\Gamma(2-\mu)W_{-\frac{3}{2}, \> \mu}\left( \frac{e^{-\frac{\pi}{2}i - z}}{2\pi nkQ^2} \right) - \sum_{m_1 = 0}^M R^{(4)}_{k,n, m_1}(z,\mu) - \sum_{m_2 = 0}^{M^\prime} R^{(4)}_{k,n, m_2}(z,\mu) \right\}.   \tag{10.11}
\end{align*}
\normalsize
Since a lemma similar to Lemma 7.1 also holds, the third and the fourth series on the right hand side in (10.11) are absolutely and uniformly convergent on every compact subset on the whole complex plane. 

Next, by the theorem of residues, (4.7) and (9.5) we have
\begin{align*}
   f_2 (z,F) + f_{2}^{-}(z,F)
  &= \int_L \frac{\zeta(s-1)}{F(s)}e^{sz}ds - \int_{\overline{L}} \frac{\zeta(s-1)}{F(s)}e^{sz}ds
   \\
   &= -2\pi i \lim_{s \to 2} (s-2) \frac{\zeta(s-1)}{F(s)}e^{zs}
   \\
   &= -2\pi i\frac{e^{2z}}{F(2)}.  \tag{10.12}
\end{align*}Finally, by (6.1) and (9.6) 
\begin{equation*}
   f_3 (z,F) + f_{3}^{-}(z,F) = 0.  \tag{10.13}
\end{equation*}Thus, for $|y| < \pi$ we have
\begin{align*}
   2\pi i(f(z,F) + f^{-} (z,F))
   &= 2\pi i\left\{ \frac{1}{2\pi i}(f_1 (z,F) + f_{1}^{-}(z,F)) - \frac{e^{2z}}{F(2)} \right\}
   \\
   &= 2\pi iB(z,F),   \tag{10.14}
\end{align*}where
\begin{equation*}
   B(z,F) = \frac{1}{2\pi i}(f_1 (z,F) + f_{1}^{-}(z,F)) - \frac{e^{2z}}{F(2)}.  \tag{10.15}
n\end{equation*}By (10.11), the function $f_1 (z,F) + f_{1}^- (z,F)$ is absolutely and uniformly convergent on every compact subset on the whole complex plane. Hence, the function $B(z,F)$ is an entire function. Since the function $f^- (z,F)$ has a meromorphic continuation for $y<\pi$ by Corollary 4.3, the function 
\[ f(z,F) = B(z,F) - f^- (z,F) \] is analytic for all $y<\pi$. Since the function $f(z,F)$ is analytic for $z \in \mathbb{H}$ and $y<\pi$, $f(z,F)$ can be analytically continued to the whole complex plane. In a similar manner, $f^- (z,F)$ can be analytically continued to the whole complex plane. Therefore, for all $z \in \mathbb{C}$ we have
\begin{equation*} 
   f(z,F) + f^{-}(z,F) =  B(z,F).  \tag{10.16}
\end{equation*}

Finally, we prove the functional equation (4.10). We recall the hypothesis that the coefficient $a_F (n)$  in the Dirichlet series of $F$ is real for all $n$. Hence, if $\rho$ is a non-trivial zero of $F$, then so is $\overline{\rho}$. For $z \in \mathbb{H}$ we have
\begin{align*}
   \overline{f(z,F)}
   &= \lim_{n \to \infty} \overline{\sum_{{ \scriptstyle \rho } \atop{\scriptstyle 0 < \text{Im} \hspace{0.01cm} \rho < T_n}} \frac{\zeta(\rho-1)}{F^\prime (\rho)}e^{\rho z}}
   \\
   &= \lim_{n \to \infty} \overline{\sum_{{ \scriptstyle \rho } \atop{\scriptstyle 0 < \text{Im} \hspace{0.01cm} \rho < T_n}} \left(\lim_{s \to \rho}\frac{s - \rho}{F(s) - F(\rho)} \zeta(s-1)e^{zs} \right)}.
\end{align*}Since $a_F (n) \in \mathbb{R}$, so $\overline{F(s)} = F(\overline{s})$ holds. Using this, we have
\begin{align*}
   f^- (\overline{z},F)
   &= \lim_{n \to \infty} \sum_{{ \scriptstyle \rho } \atop{\scriptstyle -T_n < \text{Im} \hspace{0.01cm} \rho < 0}} \frac{\zeta(\rho-1)}{F^\prime (\rho)}e^{\rho \overline{z}}
   \\
   &=  \lim_{n \to \infty} \overline{\sum_{{ \scriptstyle \rho } \atop{\scriptstyle -T_n < \text{Im} \hspace{0.01cm} \rho < 0}} \left(\lim_{s \to \rho}\frac{\overline{s} - \overline{\rho}}{\overline{F(s)} - \overline{F(\rho})} \overline{\zeta(s-1)}e^{z\overline{s}}\right)}
   \\
   &=  \lim_{n \to \infty} \overline{\sum_{{ \scriptstyle \rho } \atop{\scriptstyle -T_n < \text{Im} \hspace{0.01cm} \rho < 0}} \left(\lim_{s \to \rho}\frac{\overline{s} - \overline{\rho}}{{F(\overline{s})} - {F(\overline{\rho}})} \zeta(\overline{s}-1)e^{z\overline{s}} \right)}
   \\
   &= \lim_{n \to \infty} \overline{\sum_{{ \scriptstyle \rho } \atop{\scriptstyle 0 < \text{Im} \hspace{0.01cm} \rho < T_n}} \left(\lim_{\overline{s} \to \overline{\rho}}\frac{ s - \rho}{F(s) - F(\rho)} \zeta(s-1)e^{zs} \right)}
   \\
   &=  \lim_{n \to \infty} \overline{\sum_{{ \scriptstyle \rho } \atop{\scriptstyle 0 < \text{Im} \hspace{0.01cm} \rho < T_n}} \left(\lim_{s \to \rho}\frac{s - \rho}{F(s) - F(\rho)} \zeta(s-1)e^{zs} \right)}
   \\
   &= \overline{f(z,F)}.
\end{align*}Therefore, we have for $z \in \mathbb{H}$
\begin{equation*}
   f(z,F) = \overline{f^- (\overline{z},F)}.  \tag{10.17}
\end{equation*}Using (10.17), we have from (10.16)
\begin{align*}
   \overline{B(\overline{z},F)}
   &= \overline{f(\overline{z},F) + f^- (\overline{z},F)}
   \\
   &= \overline{f(\overline{z},F)} + \overline{f^- (\overline{z},F)}
   \\
   &=  \overline{f(\overline{z},F)} + f(z,F).
\end{align*}Using (10.17) and (10.16) again, we have for $z \in \mathbb{H}$
\begin{equation*} 
   \overline{f(\overline{z},F)} + f(z,F) = f^- (z,F) + f(z,F) = B(z,F).   \tag{10.18}
\end{equation*}Since $f(z,F), f^- (z,F)$ and $B(z,F)$ are entire functions, and (10.18) holds for all $z \in \mathbb{H}$, (10.18) holds for all $z \in \mathbb{C}$ by the analytic continuation. Therefore, the functional equation (4.10) holds for all $z \in \mathbb{C}$ and we have Theorem 4.4. $\square$  

\textbf{Acknowledgments.}The author expresses his sincere gratitude to Prof. Kohji Matsumoto and Dr. Sh\={o}ta Inoue for giving me various advice in writing this paper.  Also, the author expresses his gratitude to laboratory members, in particular, to master's student Keita Nakai giving the comment on the path of integration $\mathscr{L}$ in Section 3.

\bigskip
　　　　　　　　　　　　　　Hideto Iwata 
                                                  
　　　　　　　　　　　　　　Graduate School of Mathematics

　　　　　　　　　　　　　　 Nagoya University
                                                   
　　　　　　　　　　　　　　 Furocho, Chikusa-ku,
                                                   
　　　　　　　　　　　　　　 Nagoya, 464-8602, Japan.

　　　　　　　　　　　　　　 \small{e-mail:d18001q@math.nagoya-u.ac.jp}

\end{document}